\documentclass[11pt]{article}

\usepackage{braket,amsfonts}

\usepackage{array}

\usepackage[caption=false]{subfig}
\usepackage{pgfplots}

\usepackage{algorithmic}

\usepackage{graphicx,epstopdf}

\usepackage{subfig}
\usepackage{amsmath}
\usepackage{amssymb}
\usepackage{graphicx}
\usepackage{dcolumn}
\usepackage{mathtools}
\usepackage{amsmath,amsfonts}
\usepackage{bm}
\usepackage{amsmath}
\usepackage{amssymb}
\usepackage{algorithm2e}
\usepackage{color}
\usepackage{float}
\usepackage{setspace}
\usepackage{tabularx}
\usepackage{algorithmic}

\usepackage{multirow}

\allowdisplaybreaks

\newcommand{\p}{\partial}
\newcommand{\f}{\frac}
\newcommand{\ds}{\displaystyle}
\newcommand{\E}{ {\mathbb{E}} }
\newcommand{\ola}{\overleftarrow}

\newcommand{\be}{\begin{equation}}
\newcommand{\ee}{\end{equation}}

\usepackage{amsopn}

\usepackage{xspace}
\usepackage{bold-extra}
\usepackage[most]{tcolorbox}

\colorlet{texcscolor}{blue!50!black}
\colorlet{texemcolor}{red!70!black}
\colorlet{texpreamble}{red!70!black}
\colorlet{codebackground}{black!25!white!25}

\newtcblisting[use counter=example]{example}[2][]{%
  title={Example~\thetcbcounter: #2},#1}

\newtcbinputlisting[use counter=example]{\examplefile}[3][]{%
  title={Example~\thetcbcounter: #2},listing file={#3},#1}

\DeclareTotalTCBox{\code}{ v O{} }
{ 
  fontupper=\ttfamily\color{black},
  nobeforeafter,
  tcbox raise base,
  colback=codebackground,colframe=white,
  top=0pt,bottom=0pt,left=0mm,right=0mm,
  leftrule=0pt,rightrule=0pt,toprule=0mm,bottomrule=0mm,
  boxsep=0.5mm,
  #2}{#1}

\patchcmd\newpage{\vfil}{}{}{}
\date{}

\begin{document}

\title{A Kernel Learning Method for Backward SDE Filter}
\author{
Richard Archibald\thanks{ Computational Science and Mathematics Division, Oak Ridge National Laboratory, Oak Ridge, Tennessee.} , 
\and Feng Bao\thanks{ Department of Mathematics, Florida State University, Tallahassee, Florida, \ ({\tt bao@math.fsu.edu}).}   }
\maketitle


\begin{abstract}
In this paper, we develop a kernel learning backward SDE filter method to estimate the state of a stochastic dynamical system based on its partial noisy observations. A system of forward backward stochastic differential equations is used to propagate the state of the target dynamical model, and Bayesian inference is applied to incorporate the observational information. To characterize the dynamical model in the entire state space, we introduce a kernel learning method to learn a continuous global approximation for the conditional probability density function of the target state by using discrete approximated density values as training data.
Numerical experiments demonstrate that the kernel learning backward SDE is highly effective and highly efficient.
\end{abstract}

\vspace{2em}




\section{Introduction}

One of the key missions in data assimilation is to obtain the best estimate for the state of a stochastic dynamical system based on its observations. The mathematical tool that achieves this mission is the optimal filtering. An optimal filtering problem is usually composed of a stochastic differential equation (SDE) called the state equation, which describes the the state of the dynamical system, and an observation equation that provides partial noisy observational data. The best estimate that we want to obtain in the optimal filtering problem is the conditional expectation of the target state, which is conditioned on the observational information. When both the dynamical system and the observations are linear, the Kalman filter provides an analytical expression for the solution of the optimal filtering problem through Bayesian inference \cite{Kalman1961}. However, in most situations, we do not have the linearity condition, hence numerical methods for solving nonlinear filtering problems are needed.

The main theme of nonlinear filtering methods is to derive approximations for the conditional probability density function (PDF) of the target state, which is also called the ``filtering density''. 
An important pioneer approach to solve the nonlinear filtering problem is the Zakai filter. In the Zakai filter, the filtering density is formulated as the solution of a parabolic type stochastic partial differential equation (SPDE) called the Zakai equation \cite{zakai}. Although the Zakai filter provides a mathematical equation that analytically solves the nonlinear filtering problem, obtaining numerical solutions for SPDEs is a challenging task \cite{HU-Zakai, Gobet-Zakai}.  Especially, when the dimension of the problem is high, numerical methods for SPDEs suffer from the ``curse of dimensionality'', and 
the computational cost for solving the Zakai equation increases exponentially as the dimension of the problem increases  \cite{Zhang_Zakai, Bao_Zakai}.

The standard approach to solve the nonlinear filtering problem in practice is the Bayes filter.  Well-known Bayes filter methods include the Kalman type filters and the particle filter. The Kalman type filters \cite{ExKF, UnKF, EvensenBook} for the nonlinear filtering problem usually linearize the nonlinear systems, and then use the Kalman filter method to solve the corresponding linearized problem. The main drawback of the Kalman type filters is that when the nonlinear systems are highly nonlinear, the linearized problem does not provide a good approximation for the original problem, hence the Kalman type filters often fail \cite{Tong_EnKF}.  The central idea of the particle filter method is to use a set of sequentially generated samples (called ``particles'') to construct an empirical distribution as a predicted filtering density. 
To incorporate the observational information, the particle filter applies Bayesian inference to assign a likelihood value to each particle as its weight, and use the weighted particles to describe the updated filtering density \cite{particle-filter}. However, since the particles are generated from a stochastic dynamical system, they diffuse in long term simulations. Therefore, the particle filter has the so-called ``degeneracy problem'', i.e., when estimating the target state for several steps, only a few particles stay in high probability regions of the filtering density, and the others lie in probabilistically insignificant regions. This makes the effective particle-size decrease dramatically. To address the degeneracy of particles, a resampling procedure is introduced to re-generate particles in high probability regions \cite{APF, CT1, MCMC-PF, vanLeeuwen, MTAC2012, BaoCH_CiCP, Kang-PF}. But when solving highly nonlinear or high dimensional problems, the existing resampling techniques are either less effective or very difficulty to implement \cite{Sny-particle}.

In a recent study, we developed a backward doubly stochastic differential equation (BDSDE) approach to solve the nonlinear filtering problem, in which we use a BDSDE system to model the filtering density \cite{Bao_first, BCZ_2011, BCZ_2018, Bao_CMS}, and a ``doubly stochastic integral'' term is introduced to incorporate the observations. In this way, the BDSDE approach is similar to the Zakai filter in the sense that it also provides a mathematical equation to formulate the analytical solution of the nonlinear filtering problem \cite{PP1994, BSDE_finance}. However, in order to solve BDSDEs, a type of two-sided It\^o formula is needed to deal with the doubly stochastic integral term \cite{PP1994, BCZ_2015, BSDE_filter, BaoCC_2019}, which makes the BDSDE approach quite complicated.


In this work, we simplify the BDSDE approach and introduce a forward backward stochastic differential equations (FBSDEs) system to generate the predicted filtering density under the Bayes filter framework \textit{without using the doubly stochastic integral term.}
Then, we apply Bayesian inference to incorporate the observational information and update the filtering density. 
We call the general approach that involves backward SDEs (or backward doubly SDEs) in solving the nonlinear filtering problem as the \textit{``backward SDE filter''}. 
The numerical implementation of the backward SDE filter aims to approximate the filtering density on adaptively selected spatial points in the state space, and the adaptive spatial points are chosen as the random samples that follow the filtering density. Although the idea of using random spatial sample points is similar to the particle filter, which utilizes particles to build an empirical distribution, we also approximate \textit{filtering density values} on sample points in the backward SDE filter. Therefore, the spatial samples, together with the filtering density values on those samples, carry more information than just the particle positions in the particle filter.

On the other hand, an approximation for the entire filtering density as a continuous PDF can provide a complete  description for the state of the target dynamics. This is even more advantageous in high dimensional problems since finite spatial samples tend to be very sparse in high dimensional state spaces.

The novel methodology that we want to develop in this work is to treat the approximated filtering density values obtained by the backward SDE filter on discrete spatial sample points as ``training data samples'' and then derive a continuous approximation for the filtering density by using machine learning methods. 
In this way, the information of filtering density on scattered spatial samples is effectively combined as a smooth distribution for the target state in the entire state space.
The machine learning model that we choose in this paper is the kernel learning method \cite{Kernel_learning, Multi-Kernel_Learning}, and we name this approach the \textit{kernel learning backward SDE filter}. Since the filtering density is a probability distribution, we use Gaussian kernels to construct the filtering density function. The centers of Gaussian kernels are chosen as spatial samples with high density values, so that the kernels can effectively cover high probability regions of the state. To make the kernel learning backward SDE filter more efficient, we introduce an implicit iterative scheme to solve backward SDEs, and we apply the stochastic approximation method \cite{Stochastic_Approximation} to reduce the computational cost of simulating the conditional expectation in our iterative scheme .

We want to mention that our adaptive spatial sample points are generated through the state dynamical model -- just like the particle generation in the particle filter. To avoid sample degeneracy, we introduce a resampling procedure to resample the spatial points from the filtering density. Note that the kernel learned filtering density is a combination of Gaussian kernels, hence the resampling procedure mainly samples Gaussian variables, which can be implemented  accurately and efficiently.

The rest of this paper is organized as follows.  In Section \ref{Prelim}, we provide some preliminaries for the nonlinear filtering problem, the Bayes filter approach, and the mathematical framework of the backward SDE filter. In Section \ref{Schemes}, we introduce numerical algorithms for the kernel learning backward SDE filter. Several numerical experiments will be carried out in Section \ref{Numerics} to demonstrate the effectiveness and efficiency of the kernel learning backward SDE filter. Some concluding remarks will be given in Section \ref{Conclusion}.

\section{Preliminaries}\label{Prelim}

In this section, we provide the necessary preliminary knowledge for the kernel learning backward SDE filter. Specifically, we first provide a brief introduction to the nonlinear filtering problem. Then, we discuss the state-of-the-art approach, i.e., the Bayes filter for solving the nonlinear filtering problem. Finally, we introduce the mathematical formulation of the backward SDE filter.

\vspace{0.3em}
\noindent \textbf{The nonlinear filtering problem.} 
\vspace{0.3em}

We consider the following state-space model of the nonlinear filtering problem in a complete probability space $(\Omega, \mathcal{F}, P)$:
\begin{subequations}\label{Intro:NF}
\begin{align}
dS_t =& \ b(S_t) dt + \sigma_t dW_t, \qquad &\text{(State)} \label{NF:State} \\
dM_t =& \ h(S_t) dt + d V_t, \qquad &\text{(Observation)} \label{NF:Measurement}
\end{align}
\end{subequations}
where $S_t \in \mathbb{R}^d$ describes the state of a stochastic dynamical system driven by the nonlinear function $b: \mathbb{R}^d \rightarrow \mathbb{R}^d$, and we assume that the initial state $S_0$ follows a given distribution $p_0$. The process $W:=\{W_t\}_{t \geq 0}$ is a standard $d$-dimensional Brownian motion with the coefficient $\sigma_t \in \mathbb{R}^{d \times d}$, and the It\^o integral term brings noises that perturb the dynamical model, and we use $\mathcal{F}_t^W := \sigma(W_s, 0 \leq s \leq t)$ to denote the $\sigma$-algebra generated by $W$. The $r$-dimensional stochastic process $M_t$ \eqref{NF:Measurement} gives partial noisy observations on $S_t$ through the nonlinear observation function $h: \mathbb{R}^d \rightarrow \mathbb{R}^r$, and the observational data are also perturbed by noises generated by an $r$-dimensional Brownian motion $V$, which is independent of $W$.  We call the equation \eqref{NF:State} the ``state equation'' and the equation \eqref{NF:Measurement} the ``observation equation''.  The goal of the nonlinear filtering problem is to obtain the ``best'' estimate for the state $S$ given the observational information $\mathcal{M}_t : = \sigma(M_s, 0 \leq s \leq t)$, which is the $\sigma$-algebra generated by the observation equation. More generally, we want to determine the so-called ``optimal filter'' $\tilde{\Psi}(S_t)$ for a test function $\Psi$ representing the quantity of interest. Mathematically, the optimal filter $\tilde{\Psi}(S_t)$ is given by the conditional expectation of $\Psi(S_t)$ given $\mathcal{M}_t$, i.e.
\begin{equation}\label{Def:filter}
\tilde{\Psi}(S_t) = \E[\Psi(S_t) \big| \mathcal{M}_t].
\end{equation}

In this work, we focus on the Bayes filter approach, which is carried out by recursive Bayesian estimations. Instead of trying to estimate the optimal filter $\tilde{\Psi}$ as a conditional expectation directly, the Bayes filter aims to approximate the conditional probability density function (PDF) of the target state, which is also called the ``filtering density''. 
 
\vspace{0.3em}
\noindent \textbf{Recursive Bayesian estimations.} 
\vspace{0.3em}

In the Bayes filter, we estimate the target state on a sequence of discrete time instants $0 = t_0 < t_1 < t_2 < \cdots < t_n < \cdots t_{N_{T}} = T$ over the time interval $[0, T]$, where $N_T \in \mathbb{N}$ is the total number of time steps. The main theme of recursive Bayesian estimations is to obtain the filtering density at the time instant $t_{n+1}$, i.e. $p(S_{t_{n+1}} \big| \mathcal{M}_{t_{n+1}})$, through recursive Bayesian inferences. Then, the optimal filter $\tilde{\Psi}(S_{t_{n+1}})$ can be calculated as $\tilde{\Psi}(S_{t_{n+1}}) = \int \Psi(S_{t_{n+1}}) p(S_{t_{n+1}} \big| \mathcal{M}_{t_{n+1}}) d S_{t_{n+1}}$. The general framework of recursive Bayesian estimations is composed of two steps: a \textit{prediction step} and an \textit{update step}. 

In the prediction step, assuming that the filtering density $p(S_{t_{n}} \big| \mathcal{M}_{t_n})$ is available at the time instant $t_n$, we use the Chapman-Kolmogorov formula to propagate the dynamical model in the state equation as follows
\begin{equation}\label{CK-prediction}
p(S_{t_{n+1}} \big| \mathcal{M}_{t_n}) = \int p(S_{t_{n+1}} \big| S_{t_n})p(S_{t_{n}} \big| \mathcal{M}_{t_n}) d S_{t_n}, 
\end{equation}
where $p(S_{t_{n+1}} \big| S_{t_n})$ is the transition probability of the state equation \eqref{NF:State}, and the predicted filtering density $p(S_{t_{n+1}} \big| \mathcal{M}_{t_n})$, which is the prior distribution in the Bayesian inference, describes the state $S$ at the time instant $t_{n+1}$ before receiving the new observational data $M_{t_{n+1}}$. 

In the update step, we use the following Bayesian inference formula to to incorporate the observational data into the state estimation:
\begin{equation}\label{Bayesian-update}
p(S_{t_{n+1}} \big| \mathcal{M}_{t_{n+1}}) = \f{p(M_{t_{n+1}} \big| S_{t_{n+1}}) p(S_{t_{n+1}} \big| \mathcal{M}_{t_n}) }{p(M_{t_{n+1}} \big| \mathcal{M}_{t_n})},
\end{equation}
where $p(M_{t_{n+1}} \big| S_{t_{n+1}})$ is the likelihood function, and the denominator $p(M_{t_{n+1}} \big| \mathcal{M}_{t_n})$ is a normalization factor.

Then, by implementing the prediction step \eqref{CK-prediction} and the update step \eqref{Bayesian-update} numerically, one can develop computational methods for solving the nonlinear filtering problem.  In what follows, we introduce the backward SDE filter as the theoretical preparation for our kernel learning method.

\vspace{0.3em}

\noindent \textbf{The backward SDE filter.} 
\vspace{0.3em}

The backward SDE filter adopts the recursive Bayesian estimations framework. The central idea of the backward SDE filter is to use a system of (time-inverse) forward backward stochastic differential equations to propagate the filtering density, and we also use Bayesian inference to incorporate the observational information into the predicted filtering density. 

To proceed, we first introduce the forward backward stochastic differential equations (FBSDEs) corresponding to the nonlinear filtering problem \eqref{Intro:NF}, and we consider the following FBSDEs system
\begin{equation}\label{FBSDE-forward}
\begin{aligned}
S_t =& S_0+ \int_{0}^t b(S_s) ds + \int_{0}^t \sigma_s dW_s, \\
\tilde{Y}_0 =& \tilde{Y}_t - \int_{0}^{t} \tilde{Z}_s  dW_s, \hspace{4em} \tilde{Y}_t = \Psi(S_t),
\end{aligned}
\end{equation}
where the first equation coincides the state equation in the nonlinear filtering problem, which is a standard forward SDE, and the second equation is a backward SDE. The solution of the above FBSDEs system is the pair $(\tilde{Y}, \tilde{Z})$, which is adapted to the Brownian motion $W$, i.e., $\tilde{Y}_t, \tilde{Z}_t  \in \mathcal{F}_{t}^{W}$, and the solution $\tilde{Z}$ is the martingale representation of $\tilde{Y}$ with respect to $W$ \cite{Peng_ICM}. Note that the side condition of the backward SDE, i.e., $\tilde{Y}_t = \Psi(S_t)$, is given at the time instant $t$, and the solution pair $(\tilde{Y}, \tilde{Z})$ propagates backwards from $t$ to $0$. 

For a fixed initial state $S_0 = x \in \mathbb{R}^d$, we take the conditional expectation $\E[\cdot \big| S_0 = x]$ on both sides of the backward SDE in \eqref{FBSDE-forward} and obtain $\tilde{Y}_0(x) = \E[\Psi(S_t) \big| S_0=x]$, which is a simplified version of the Feynman-Kac formula. Here, we note that the value of $\tilde{Y}_0$ is determined by the value of the state $S_0$. 
In addition, the solution $\tilde{Y}$ of \eqref{FBSDE-forward} is equivalent to the solution of the Kolmogorov backward equation. In other words, for the following \textit{backward} parabolic type partial differential equation (PDE),
\begin{equation}\label{Kolmogorov}
- \f{d u_s}{ds} = \sum_{i=1}^d b_i \f{\p u_s}{\p x_i} + \f{1}{2} \sum_{i,j=1}^{d} (\sigma_s \sigma_s^T)_{i, j} \f{\p^2 u_s}{\p x_i \p x_j},  \qquad u_t(x) = \Psi(x),
\end{equation}
where $b_i$ is the $i$-th component of the vector function $b$ in the state equation \eqref{NF:State}, we have $\tilde{Y}_0(S_0=x) = u_0(x)$ \cite{Pardoux1991}.

In the nonlinear filtering problem, we need to propagate the filtering density forward. The PDE that propagates the PDF of the state $S$ driven by \eqref{NF:State} is the following Fokker-Planck equation:
\begin{equation}\label{FP-equation}
\f{d p_t}{dt} = - \sum_{i=1}^d \f{\p b_i \ p_t}{\p x_i} + \f{1}{2} \sum_{i,j=1}^{d} (\sigma_t \sigma_t^T)_{i, j} \f{\p^2 p_t}{\p x_i \p x_j},
\end{equation}
where the initial condition $p_0$ is the distribution of the state $S_0$. We can see that \eqref{FP-equation} is the adjoint equation of the Kolmogorov backward equation, hence the Fokker-Planck equation is also called the Kolmogorov forward equation.

Following the analysis in \cite{Pardoux1991} that establishes the equivalence between the  FBSDEs \eqref{FBSDE-forward} and the Kolmogorov backward equation \eqref{Kolmogorov}, one can derive that the solution $p_t$ of the above Fokker-Planck equation is equivalent to the solution $Y_t$ of the following FBSDEs:
\begin{subequations}\label{FBDSDEs:Adjoint}
\begin{align}
X_0 &= X_t - \int_{0}^{t} b(X_s) ds + \int_{0}^{t} \sigma_s d\ola{W}_s,    \label{FBDSDEs:Adjoint:a} \\
Y_t &= Y_0 - \int_{0}^{t} \sum_{i=1}^{d}\f{\p b_i}{\p x_i}(X_s) Y_s  ds  - \int_{0}^{t} Z_s d\ola{W}_s, \qquad  Y_0= p_0, \label{FBDSDEs:Adjoint:b}
\end{align}
\end{subequations}
where the integral $ \int_{0}^{t} \cdot d\ola{W}_s$ is a \textit{backward} It\^o integral, which is an It\^o integral integrated backwards  \cite{PP1994}. We can see that the first equation in \eqref{FBDSDEs:Adjoint} is a ``backward SDE'' since it propagates backwards from $t$ to $0$, and the second equation is a ``forward SDE'' since its side condition is given at the time instant $0$. On the other hand, the propagation direction of \eqref{FBDSDEs:Adjoint:a} is the same as the integration direction of the backward It\^o integral. At the same time, the equation \eqref{FBDSDEs:Adjoint:b} propagates forwards given the side condition $Y_0=p_0$ with a backward It\^o integral. In this way, the equations in \eqref{FBDSDEs:Adjoint} compose a \textit{time-inverse} FBSDEs system. 

Due to the equivalence $Y_t = p_t$, we know that the solution $Y_t$ also propagates the PDF of the state $S_{t}$ forward.  In this connection, we use the time-inverse FBSDEs system \eqref{FBDSDEs:Adjoint} to predict the filtering density in the nonlinear filtering problem.
Specifically, we assume that we have the filtering density $p(S_{t_n} \big| \mathcal{M}_{t_n})$ at the time step $t_n$.  
By solving the following time-inverse FBSDEs system
\begin{subequations}\label{FBDSDEs:step}
\begin{align}
X_{t_n} &= X_{t_{n+1}} - \int_{t_n}^{t_{n+1}} b(X_s) ds + \int_{t_n}^{t_{n+1}} \sigma_s d\ola{W}_s,    \label{FBDSDEs:step:a} \\
Y^{M_{t_n}}_{t_{n+1}} &= Y^{M_{t_n}}_{t_n} - \int_{t_n}^{t_{n+1}} \sum_{i=1}^{d}\f{\p b_i}{\p x_i}(X_s) Y^{M_{t_n}}_s  ds  - \int_{t_n}^{t_{n+1}} Z^{M_{t_n}}_s d\ola{W}_s,  \label{FBDSDEs:step:b}\\
X_{t_{n+1}} & = \ S_{t_{n+1}}, \qquad Y^{M_{t_n}}_{t_n} = p(S_{t_n}\big|\mathcal{M}_{t_n}), \nonumber
\end{align}
\end{subequations}
we obtain the solution $Y^{M_{t_n}}_{t_{n+1}}$, which is the predicted filtering density, i.e., $Y^{M_{t_n}}_{t_{n+1}} = p(S_{t_{n+1}}\big| \mathcal{M}_{t_n})$.
In other words, the time-inverse FBSDEs system \eqref{FBDSDEs:step} provides a mechanism to carry out the the Chapman-Kolmogorov formula \eqref{CK-prediction} in the prediction step of the Bayes filter approach.  Here, we use the superscript $M_{t_n}$ in $Y^{M_{t_n}}_{t_{n+1}}$ to emphasize that the solution $Y$ depends on the observational information $\mathcal{M}_{t_n}$. 

Then, we substitute the prior distribution in the Bayesian inference \eqref{Bayesian-update} by the solution $Y^{M_{t_n}}_{t_{n+1}}$ obtained in \eqref{FBDSDEs:step} to get the posterior distribution $p(S_{t_{n+1}}\big| \mathcal{M}_{t_{n+1}})$, i.e.
\begin{equation}\label{Bayesian-Y}
p(S_{t_{n+1}} \big| \mathcal{M}_{t_{n+1}}) = \f{p(M_{t_{n+1}} \big| S_{t_{n+1}}) Y^{M_{t_n}}_{t_{n+1}} }{p(M_{t_{n+1}} \big| \mathcal{M}_{t_n})},
\end{equation}
 which is the updated filtering density at time step $t_{n+1}$ that will be used for the next recursive stage. 

From the above discussion, we can see that the central idea of the backward SDE filter is to use the time-inverse FBSDEs system to predict the filtering density, and then use Bayesian inference to update the predicted filtering density. As a method that carries out recursive Bayesian estimations, the backward SDE filter is also composed of a prediction step and an update step. In most practical situations, FBSDEs are not explicitly solvable. Therefore, numerical solutions for FBSDEs are needed.
In the following section, we introduce a numerical algorithm to implement the above backward SDE filter framework and develop our efficient kernel learning method for the backward SDE filter.


\section{Numerical algorithms for kernel learning backward SDE filter}\label{Schemes}


We first provide numerical schemes to solve the time-inverse FBSDEs system \eqref{FBDSDEs:step} and then give Bayesian inference based on the numerical solution of the FBSDEs system in Subsection \ref{Num:BSDE}.  In order to approximate the entire filtering density that carries the information contained in the state dynamical model and the observational data, in Subsection \ref{KL-BSDEF} we introduce a kernel learning method to ``learn'' the filtering density from discrete density values. In Subsection \ref{Summary} we summarize our kernel learning backward SDE filter.

\subsection{Numerical schemes for time-inverse FBSDEs}\label{Num:BSDE}


Since the equation \eqref{FBDSDEs:step:a} is essentially an SDE with inverse propagation direction, we apply the Euler-Maruyama scheme \cite{Kloeden1992} and get
\begin{equation*}
X_{t_n} = X_{t_{n+1}} - b(X_{t_{n+1}}) \Delta t_n + \sigma_{t_{n+1}} \Delta W_{t_n} + R_n^X,
\end{equation*}
where $\Delta t_n := t_{n+1} - t_n$, $\Delta W_{t_n} := W_{t_{n+1}} - W_{t_n}$, and $R_n^X := b(X_{t_{n+1}}) \Delta t_n - \int_{t_n}^{t_{n+1}} b(X_s) ds$ is the approximation error. By dropping the error term $R_n^X$, we obtain the following numerical scheme for \eqref{FBDSDEs:step:a}:
\begin{equation}\label{Scheme:SDE}
X_{n} = X_{n+1} - b(X_{n+1}) \Delta t_n + \sigma_{t_{n+1}} \Delta W_{t_n},
\end{equation}
where $X_{n}$ is the numerical approximation for $X_{t_n}$, and $X_{n+1}$ is a representation for the state variable $S_{t_{n+1}}$.

To solve the backward SDE \eqref{FBDSDEs:step:b}, we take conditional expectation $\E_{n+1}^X[\cdot] := \E[\cdot \big| X_{t_{n+1}}, \mathcal{M}_{t_n}]$ on both sides of the equation and obtain
\begin{equation}\label{Exp:BSDE}
\E_{n+1}^X[ Y^{M_{t_n}}_{t_{n+1}}] = \E_{n+1}^X[ Y^{M_{t_n}}_{t_n}] - \int_{t_n}^{t_{n+1}} \E_{n+1}^X\Big[ \sum_{i=1}^{d}\f{\p b_i}{\p x_i}(X_s) Y^{M_{t_n}}_s \Big] ds,
\end{equation}
where the backward It\^o integral $\int_{t_n}^{t_{n+1}} Z_s d\ola{W}_s$ is eliminated due to the martingale property of It\^o integrals. For the left hand side of the above equation, we have $Y^{M_{t_n}}_{t_{n+1}} = \E_{n+1}^X[ Y^{M_{t_n}}_{t_{n+1}}]$ since  $Y^{M_{t_n}}_{t_{n+1}}$ is adapted to $X_{t_{n+1}}$ and is $\mathcal{M}_{t_n}$ measurable. 
In this work, we approximate the deterministic integral on the right hand side of \eqref{Exp:BSDE} by using the right-point formula and get
\begin{equation}\label{Implicit:BSDE}
Y^{M_{t_n}}_{t_{n+1}} = \E_{n+1}^X[ Y^{M_{t_n}}_{t_n}] -  \Delta t_n \cdot \sum_{i=1}^{d}\f{\p b_i}{\p x_i}(X_{t_{n+1}}) Y^{M_{t_n}}_{t_{n+1}} + R_n^Y,
\end{equation}
where $R_n^Y := \Delta t_n \cdot \sum_{i=1}^{d}\f{\p b_i}{\p x_i}(X_{t_{n+1}}) Y^{M_{t_n}}_{t_{n+1}}  - \int_{t_n}^{t_{n+1}} \E_{n+1}^X\Big[ \sum_{i=1}^{d}\f{\p b_i}{\p x_i}(X_s) Y^{M_{t_n}}_s \Big] ds$ is the approximation error for the integral, and we have used the fact $\sum_{i=1}^{d}\f{\p b_i}{\p x_i}(X_{t_{n+1}}) Y^{M_{t_n}}_{t_{n+1}} = \E_{n+1}^X\big[\sum_{i=1}^{d}\f{\p b_i}{\p x_i}(X_{t_{n+1}}) Y^{M_{t_n}}_{t_{n+1}}\big]$. 

Then, we drop the approximation error term $R_n^Y$ in \eqref{Implicit:BSDE} and obtain the following approximation scheme for $Y^{M_{t_n}}_{t_{n+1}}$:
\begin{equation}\label{Approx:Y}
Y^{M_{t_n}}_{n+1} = \E_{n+1}^X[ Y^{M_{t_n}}_{n}] -  \Delta t_n \cdot \sum_{i=1}^{d}\f{\p b_i}{\p x_i}(X_{t_{n+1}}) Y^{M_{t_n}}_{n+1},
\end{equation}
where $Y^{M_{t_n}}_{n+1}$ is the approximated solution and $Y^{M_{t_n}}_{n}$ is an approximation of the filtering density $p(S_{t_n}\big|\mathcal{M}_{t_n})$ that we obtained in the previous recursive step.
We can see that the above approximation scheme is an implicit scheme. In order to calculate $Y^{M_{t_n}}_{n+1}$, we introduce the following fixed-point iteration procedure:
\begin{equation}\label{Scheme:Y:Semi}
Y_{n+1}^{M_{t_n}, l+1} =  \E_{n+1}^X[Y_{n}^{M_{t_n}}] -  \Delta t_n \cdot \sum_{i=1}^{d}\f{\p b_i}{\p x_i}(X_{t_{n+1}}) Y^{M_{t_n}, l}_{n+1}, \quad l = 0, 1, 2, \cdots, L-1,
\end{equation}
and we let the approximated solution be $Y_{n+1}^{M_{t_n}} = Y_{n+1}^{M_{t_n}, L}$, where $L$ is a number that satisfies certain stopping criteria for the iterations, and we let the initial guess for the solution $Y_{n+1}^{M_{t_n}}$ be $Y_{n+1}^{M_{t_n}, 0} = Y_{n}^{M_{t_n}}$

From the equivalence $Y_{t_{n+1}}^{M_{t_n}} = p(S_{t_{n+1}}\big| \mathcal{M}_{t_n})$, the approximated solution $Y_{n+1}^{M_{t_n}}$ gives an approximation for the predicted filtering density. Hence the iterative scheme \eqref{Scheme:Y:Semi} accomplishes the prediction step in Bayesian estimation. To incorporate the observational information, we carry out the update step through Bayesian inference as follows
\begin{equation}\label{Bayesian-Semi}
\tilde{p}(S_{t_{n+1}} \big| \mathcal{M}_{t_{n+1}}) = \f{p(M_{t_{n+1}} \big| S_{t_{n+1}}) Y_{n+1}^{M_{t_n}} }{C},
\end{equation}
where $C$ is a normalization factor, and the prior distribution is replaced by the approximation of the predicted filtering density $Y_{n+1}^{M_{t_n}}$.  As a result, we obtain the approximated filtering density $\tilde{p}(S_{t_{n+1}} \big| \mathcal{M}_{t_{n+1}}) $ as desired in the backward SDE filter. 

The numerical schemes \eqref{Scheme:Y:Semi}-\eqref{Bayesian-Semi} compose a general computational framework for the backward SDE filter, which provides a recursive prediction-update mechanism that formulates the temporal propagation of the filtering density. On the other hand, the filtering density is a function that connects state positions to the probability density values at those positions. Therefore, spatial dimension approximation for the filtering density function with respect to the state variable is needed. 
In the following subsection, we introduce a kernel learning method to generate a continuous global approximation for the filtering density over the state space. 

%




\subsection{Efficient kernel learning in the backward SDE filter}\label{KL-BSDEF}
To derive an efficient kernel learning method that approximates the filtering density in the backward SDE filter, we re-consider the temporal prediction-update schemes \eqref{Scheme:Y:Semi}-\eqref{Bayesian-Semi}. We can see that in the scheme \eqref{Scheme:Y:Semi}, the conditional expectation $\E_{n+1}^X[Y_{n}^{M_{t_n}}]$ needs to be evaluated in order to calculate the predicted filtering density $Y_{n+1}^{M_{t_n}}$.  In what follows, we first discuss our approximation method for $\E_{n+1}^x[Y_{n}^{M_{t_n}}]$ given that $X_{t_{n+1}} = x$. 

\vspace{0.5em}
\noindent \textbf{Approximating the conditional expectation} 
\vspace{0.5em}

In most practical nonlinear filtering problems, the target state is a high dimensional variable, and Monte Carlo simulation is usually applied to evaluate high dimensional expectations. Specifically, for a given point $x \in \mathbb{R}^d$ in the state space, we let $X_{n+1} = x$ in the scheme \eqref{Scheme:SDE}. Then, the conditional expectation in \eqref{Scheme:Y:Semi} is approximated by
\begin{equation}\label{MC-approximation}
\E_{n+1}^x[Y_{n}^{M_{t_n}}] \approx  \hat{\E}_{n+1}^x[Y_{n}^{M_{t_n}}] := \f{\sum_{m=1}^M Y_{n}^{M_{t_n}} (\tilde{X}_n^{x, m})}{M},
\end{equation}
where $M$ is the total number of Monte Carlo samples that we use to approximate the conditional expectation. The random variable $\tilde{X}_n^{x, m}$ in the above approximation is a Monte Carlo sample simulated by using the scheme \eqref{Scheme:SDE} as follows:
\begin{equation}\label{X:sample-m}
\tilde{X}_{n}^{x, m} = x - b(x) \Delta t_n + \sigma_{t_{n+1}} \sqrt{\Delta t_n} \omega_m,
\end{equation}
where $\omega_m$ is a sample drawn from the $d$-dimensional standard Gaussian distribution, and $\sqrt{\Delta t_n} \omega_m$ is the $m$-th realization of $\Delta W_{t_n}$. 
However, when the dimension of the state variable is high, the number $M$ of Monte Carlo samples needs to be very large, hence evaluating the conditional expectation $\E_{n+1}^x[Y_{n}^{M_{t_n}}]$  by using the Monte Carlo simulation \eqref{MC-approximation} is a computationally expensive task.   

In this work, inspired by the stochastic approximation method and its application in stochastic optimization \cite{Stochastic_Approximation, SGD_NIPS2007}, we treat the large number of Monte Carlo samples that we use to approximate the conditional expectation $\E_{n+1}^x[Y_{n}^{M_{t_n}}]$ as a ``large data set''. Then, we adopt the methodology of stochastic approximation and use a single-sample (or a small batch of samples) to represent conditional expectations. Specifically, in each fixed-point iteration step \eqref{Scheme:Y:Semi}, instead of using the fully-calculated Monte Carlo simulation \eqref{MC-approximation} to compute the conditional expectation, one may use one realization of the simulated sample $\tilde{X}_{n}^{x, l}$ to represent the entire set of Monte Carlo samples $\{\tilde{X}_n^{x, m}\}_{m=1}^M$ at each iteration step, where $\tilde{X}_{n}^{x, l}$ is also simulated through \eqref{X:sample-m} indexed by the iteration step $l$. In this way, the fixed-point iteration scheme for the approximated solution $Y_{n+1}^{M_{t_n}}$ at the spatial point $X_{n+1} = x$ can be carried out as follows:
\begin{equation}\label{Scheme:Y:single}
Y_{n+1}^{M_{t_n}, l+1}(x) =  \tilde{\E}_{n+1}^{x, l}[Y_{n}^{M_{t_n}}]  -  \Delta t_n \cdot \sum_{i=1}^{d}\f{\p b_i}{\p x_i}(x) Y^{M_{t_n}, l}_{n+1}, \quad l = 0, 1, 2, \cdots, L-1,
\end{equation}
where the conditional expectation $\E_{n+1}^x[Y_{n}^{M_{t_n}}]$ is represented by a single-sample of $Y_{n}^{M_{t_n}}$ corresponding to $\tilde{X}_{n}^{x, l}$, and we have
\begin{equation}\label{singel-sample}
\tilde{\E}_{n+1}^{x, l}[Y_{n}^{M_{t_n}}] = Y_{n}^{M_{t_n}}(\tilde{X}_{n}^{x, l}).
\end{equation}
As a result, in each fixed-point iteration step, we only need to generate one sample of $X_n$ and evaluate the function value of the previous filtering density $Y_{n}^{M_{t_n}}$ at the spatial point $\tilde{X}_{n}^{x, l}$. In this way, we transfer the cost of simulating a large number of  Monte Carlo samples to carrying out fixed-point iterations. Although the single-sample representation does not provide accurate approximation for the conditional expectation, every simulated sample $\tilde{X}_{n}^{x, l}$ is effectively used to improve the estimate of the desired predicted filtering density $Y_{n+1}^{M_{t_n}}(x)$, which makes the overall fixed-point iteration procedure more efficient.

\vspace{0.2em}

In order to use the simulated samples more effectively and to make the approximation for the conditional expectation more accurate, we modify the single-sample representation \eqref{singel-sample} and use the batch of samples $\{Y_{n}^{M_{t_n}}(\tilde{X}_{n}^{x, l})\}_{l=0}^{L-1}$ to approximate the expectation. Precisely, we let the approximated expectation in \eqref{Scheme:Y:single} at each iteration step be
\begin{equation}\label{batch-sample}
\tilde{\E}_{n+1}^{x,l}[Y_{n}^{M_{t_n}}] = \f{\sum_{i=1}^{l}Y_{n}^{M_{t_n}}(\tilde{X}_{n}^{x, i})}{l}.
\end{equation}
In this way, all the samples previously generated are used to evaluate the expectation. Note that the expectation $\E_{n+1}^{x}[Y_{n}^{M_{t_n}}]$ in the fixed-point iteration is independent of the estimation for the solution $Y_{n+1}^{M_{t_n}}$. Hence using more samples to approximate the expectation at each iteration step only makes the iterative scheme more accurate. 

\vspace{0.5em}
\noindent \textbf{Approximating the filtering density on random spatial sample points} 
\vspace{0.5em}

By using the iterative scheme \eqref{Scheme:Y:single}, we can calculate the predicted filtering density $Y_{n+1}^{M_{t_n}}$ on the state point $x$. Then, through the Bayesian inference scheme \eqref{Bayesian-Semi}, we can obtain an approximation for the updated filtering density $\tilde{p}(S_{t_{n+1}}=x \big| \mathcal{M}_{t_{n+1}})$ at the time step $t_{n+1}$.

In order to provide a complete description for the filtering density, we need to approximate the conditional PDF of the target state as a mapping from the state variable to PDF values. 
Standard function approximation methods use tensor-product grid points, on which we approximate function values, and then use polynomial interpolation to construct an interpolatory approximation for the entire function. If the dimension of the problem is moderately high,  sparse-grid methods are often adopted as efficient alternatives to tensor-product grid interpolations.  However, even advanced adaptive sparse-grid methods suffer from the ``curse of dimensionality'' problem.  When the dimension of the problem is higher, i.e., $d \geq 10$, the cost of implementing sparse-grid approximation becomes extremely high. In many practical nonlinear filtering problems, the state dimensions are very high. Hence, applying traditional grid-based function approximation methods for filtering density is infeasible in solving high dimensional real-world problems. 

An advantage of the backward SDE filter is that it allows us to approximate the filtering density on \textit{any} point $x$ in the state space. In this way, we don't have to solve the problem on pre-determined meshes. Instead, in this work we use randomly generated state samples as our spatial points, and we generate spatial sample points so that they adaptively follow the conditional distribution of the target state. Specifically, assuming that we have a set of spatial points $\{x^n_i\}_{i=1}^N$ that follow the previous approximated filtering density $\tilde{p}(S_{t_n} \big| \mathcal{M}_{t_n})$, where $N$ is the total number of spatial points, we propagate those spatial samples through the state dynamics \eqref{NF:State} by using the following Euler-Maruyama scheme:
\begin{equation}\label{SDE:sample}
\tilde{x}_i^{n+1} = x^n_i + b(x^n_i) \Delta t_n + \sigma_{t_n} \sqrt{\Delta t_n} \omega^i,  \qquad i = 1, 2, \cdots, N,
\end{equation}
where $\{\omega^i\}_{i=1}^N$ is a sequence of i.i.d. standard $d$-dimensional Gaussian random variables.
As a result, the sample set $\{\tilde{x}_i^{n+1}\}_{i=1}^N$ forms an empirical distribution for the prior distribution  $p(S_{t_{n+1}} \big| \mathcal{M}_{t_n})$. Then, we solve the time-inverse FBSDEs \eqref{FBDSDEs:step} on those spatial sample points through the scheme  \eqref{Scheme:Y:single} to get $\big\{Y_{n+1}^{M_{t_n}}(\tilde{x}_i^{n+1})\big\}_{i=1}^N$. When the new observational data $M_{t_{n+1}}$ is available, we use Bayesian inference to update the predicted filtering density and get $\big\{\tilde{p}(S_{t_{n+1}} = \tilde{x}_i^{n+1} \big| \mathcal{M}_{t_{n+1}}) \big\}_{i=1}^N$ as our approximations for the updated filtering density. 
In this way, the approximations on the scattered sample points $\{\tilde{x}_i^{n+1}\}_{i=1}^N$ provide a partial description for the desired filtering density. 

On the other hand, in our iterative scheme \eqref{Scheme:Y:single} we need the function value of $Y_{n}^{M_{t_n}}$ on the spatial point $\tilde{X}_{n}^{\tilde{x}_i^{n+1}, l}$, which is calculated from the scheme \eqref{X:sample-m} by choosing $x = \tilde{x}_i^{n+1}$. Apparently, $\tilde{X}_{n}^{\tilde{x}_i^{n+1}, l}$ is unlikely to be one of the existing sample points, on which we have approximated the function values for $Y_{n}^{M_{t_n}}$.  Therefore, we need to derive an approximation for the filtering density over the entire state space. Since the filtering density is approximated on random spatial points, meshfree methods are needed. Although traditional meshfree interpolation methods, such like the moving least square method and the radial basis function interpolation method, could compute interpolatory approximation based on density function values at nearby spatial points, calculating the filtering density at each point separately is computationally expensive. Especially, when the dimension of the problem is high, we have to approximate the filtering density on a very large number of spatial points, which makes local approximation methods very difficult to implement.

In what follows, we introduce a kernel learning method to ``learn'' a global approximation for the entire filtering density .


\vspace{0.5em}
\noindent \textbf{Kernel learning for the filtering density} 
\vspace{0.5em}

The kernel machine utilizes the combination of a set of pre-chosen kernels to represent a model \cite{Kernel_learning}.  
For a target model in the form of a function $F$, the kernel learning method approximates the function as $F(x) \approx g\big(\sum_{k=1}^{K} \alpha_k \phi_k(x) + \beta \big)$, $x \in \mathcal{R}^d$, where $\{\phi_k\}_{k=1}^K$ is a set of $K$ kernels, $g$ is an optional nonlinearity, $\{\alpha_k\}_{k=1}^K$ are weights of kernels, and $\beta$ is a bias parameter.
In this work, we apply the kernel learning method to approximate the filtering density under the backward SDE filter framework. 
Since our target function in the nonlinear filtering problem is a probability distribution, which is nonnegative and often bell-shaped, we drop the nonlinear function $g$ and the bias $\beta$ in the kernel learning model, and we choose Gaussian type functions as our kernels. Specifically, at the time step $t_{n+1}$, we use the following kernel learning scheme to formulate a global approximation for the filtering density $p(S_{t_{n+1}}\big| \mathcal{M}_{t_{n+1}})$:
\begin{equation}\label{KL-Pn}
p_{n+1}(x) := \sum_{k=1}^{K} \alpha^{n+1}_k \phi^{n+1}_k(x), \qquad x \in \mathcal{R}^d,
\end{equation}
where $p_{n+1}$ is the kernel learned filtering density at the time step $t_{n+1}$, the Gaussian type kernel is chosen as $\ds \phi_k^{n+1} (x) = \exp\big(- (\hat{x}_k^{n+1} - x)^2 /(\lambda^{n+1}_k )^2  \big)$ with center $\hat{x}^{n+1}_k$ and covariance $\lambda^{n+1}_k$,  and the coefficient $\alpha^{n+1}_k > 0$ is the weight of the $k$-th Gaussian kernel $\phi_k^{n+1}$. 
We can see from the scheme \eqref{KL-Pn} that the features of the kernel learned density $p_{n+1}$ depend on the choice of kernel centers $\{\hat{x}^{n+1}_k\}_{k=1}^{K}$ and the parameters $(\alpha^{n+1}, \lambda^{n+1})$, where $\alpha^{n+1} : = (\alpha^{n+1}_1, \alpha^{n+1}_2, \cdots, \alpha^{n+1}_K)^T$ and  $\lambda^{n+1} : = (\lambda^{n+1}_1, \lambda^{n+1}_2, \cdots, \lambda^{n+1}_K)^T$ denote all the weights and covariances of kernels.

Since the state sample points that we generate through the scheme \eqref{SDE:sample} adaptively follow the filtering density, which provide a good representation for the target PDF, we choose the kernel centers as a subset of the state samples $\{\tilde{x}_i^{n+1}\}_{i=1}^N$. Note that we have the value of the approximated filtering density $\tilde{p}(S_{t_{n+1}} \big| \mathcal{M}_{t_{n+1}})$ on each sample point $S_{t_{n+1}} = \tilde{x}_i^{n+1}$. Therefore, we can choose state samples with high density values as kernel centers. To avoid using too many samples in the mode of the distribution and to capture more features of the filtering density, we use importance sampling to choose $\{\hat{x}^{n+1}_k\}_{k=1}^K$ \cite{particle-filter} instead of only using samples with highest density values.

To determine the parameters of the kernels in kernel learning, we use the approximated filtering density values $\big\{\tilde{p}(S_{t_{n+1}} = \tilde{x}_i^{n+1} \big| \mathcal{M}_{t_{n+1}}) \big\}_{i=1}^N$ as simulated ``training data'', and we aim to find kernel parameters so that the kernel learned filtering density $p_{n+1}$ matches the training data.
To this end, we implement stochastic gradient descent optimization to determine the parameters $\alpha^{n+1}$ and $\lambda^{n+1}$. Specifically, we define the loss function to be minimized as
\begin{equation*}
\begin{aligned}
F^{n+1}_{\alpha, \lambda}: = & \E\Big[\big(p_{n+1}(S_{t_{n+1}}) - \tilde{p}(S_{t_{n+1}} \big| \mathcal{M}_{t_{n+1}}) \big)^2 \Big] \\
= & \E\Big[\big(\sum_{k=1}^{K} \alpha^{n+1}_k \phi^{n+1}_k(S_{t_{n+1}}) - \tilde{p}(S_{t_{n+1}} \big| \mathcal{M}_{t_{n+1}}) \big)^2 \Big] .
\end{aligned}
\end{equation*}
Since the state $S_{t_{n+1}}$ is represented by spatial samples $\{\tilde{x}_i^{n+1}\}_{i=1}^N$, the fully calculated Monte Carlo simulation for the above loss function is given as 
$$F^{n+1}_{\alpha, \lambda} \approx \f{1}{N} \sum_{i=1}^{N}\Big( p_{n+1}(S_{t_{n+1}} = \tilde{x}_i^{n+1}) - \tilde{p}(S_{t_{n+1}} = \tilde{x}_i^{n+1} \big| \mathcal{M}_{t_{n+1}})  \Big)^2. $$
Then, for pre-chosen initial estimates $\alpha^{n+1}(0)$ and $\lambda^{n+1}(0)$, we carry out the following stochastic gradient descent iteration to search for the parameters $\alpha^{n+1}$ and $\lambda^{n+1}$:
\begin{subequations}\label{Kernel-Learning}
\begin{align}
\alpha^{n+1}(j+1) = \alpha^{n+1}(j) - \rho_{\alpha}^j \nabla_{\alpha} F^{n+1}_{\alpha, \lambda}\big|_{S_{t_{n+1}} = \bar{x}_j}, \qquad j = 0, 1, 2, \cdots, J-1, \\
\lambda^{n+1}(j+1) = \lambda^{n+1}(j) - \rho_{\lambda}^j \nabla_{\lambda} F^{n+1}_{\alpha, \lambda}\big|_{S_{t_{n+1}} = \bar{x}_j}, \qquad j = 0, 1, 2, \cdots, J-1,
\end{align}
\end{subequations}
where $\rho_{\alpha}^j$ and $\rho_{\lambda}^j$ are learning rates for the parameters $\alpha$ and $\lambda$, respectively, and $J$ is the total number of iterations corresponding to a stopping criteria.  The gradients $ \nabla_{\alpha}  F^{n+1}_{\alpha, \lambda}\big|_{S_{t_{n+1}} = \bar{x}_j}$ and $ \nabla_{\lambda} F^{n+1}_{\alpha, \lambda}\big|_{S_{t_{n+1}} = \bar{x}_j}$ of the cost function $F^{n+1}_{\alpha, \lambda}$ are single-sample representations of the gradients $ \nabla_{\alpha}  F^{n+1}_{\alpha, \lambda}$ and $ \nabla_{\lambda} F^{n+1}_{\alpha, \lambda}$ by choosing a specific state sample $\bar{x}_j$ for the state variable $S_{t_{n+1}}$, where the sample $\bar{x}_j$ is picked among the sample set $\{\tilde{x}_i^{n+1}\}_{i=1}^N$. As a result, we improve the estimates for $\alpha^{n+1}$ and $\lambda^{n+1}$ gradually by comparing the kernel learned filtering density $p_{n+1}$ with approximated filtering density values on samples $\{\tilde{x}_i^{n+1}\}_{i=1}^N$. 

In order to use the approximated filtering density values more effectively, instead of picking samples uniformly from the sample set, we use importance sampling to choose samples according to their density values. In this way, it is more likely to consider higher filtering density values in the optimization procedure, which makes the stochastic gradient descent procedure more efficient. In this work, we let the covariance matrices for Gaussian type kernels be diagonal to reduce the dimension of optimization, and note that the scattered kernel centers also provide covariant features of the target filtering density.

\vspace{0.5em}
\noindent \textbf{Resampling random spatial points} 
\vspace{0.5em}

In the nonlinear filtering problem, the state equation is a diffusion process. Thus the random spatial samples propagated through the scheme \eqref{SDE:sample} diffuse after several estimation steps. 
As a result, fewer and fewer spatial samples will remain in high probability regions of the filtering density as we estimate the target state step-by-step. To rejuvenate the spatial samples and to make them better represent the filtering density, we carry out a resampling procedure in the kernel learning backward SDE filter.

To be specific, we use the kernel learned updated filtering density $p_{n+1}$ to generate a set of new spatial samples, denoted by $\{x_i^{n+1}\}_{i=1}^N$, to replace the samples $\{\tilde{x}_i^{n+1}\}_{i=1}^N$, which follow the predicted filtering density. We want to point out that the kernel learned filtering density provides the conditional PDF for the target state over the entire state space. Therefore, our resampling procedure also allows us to consider probabilistically insignificant regions. On the other hand, the filtering density $p_{n+1}$ is essentially a combination of Gaussian kernels. Hence drawing samples from $p_{n+1}$ is very efficient. For example, to generate the sample $x_i^{n+1}$, we first use importance sampling to pick a kernel $\phi_k^{n+1}$ based on weights of kernels.  Since $\phi_k^{n+1}$ is a Gaussian kernel, we can simply draw the sample $x_i^{n+1}$ from the Gaussian distribution $N(\hat{x}_k^{n+1}, \lambda_k^{n+1})$.

\subsection{Summary of the algorithm}\label{Summary}

In Table \ref{Algorithm}, we summarize our kernel learning backward SDE filter as a pseudo-algorithm.
\vspace{0.5em}

\begin{table}[h!]\caption{}\label{Algorithm}
\centering
\begin{tabular} {p{0.9\textwidth}}
\hline\noalign{\smallskip}
{\bf Algorithm}: {\em Kernel learning backward SDE filter}\\
\noalign
{\smallskip}\hline
\noalign{\smallskip}
\vspace{-0.1cm}
\begin{spacing}{1.1}
\begin{algorithmic}\label{algorithm}
\item[Initialize] the spatial sample cloud $\{x_i^{0}\}_{i=1}^N \sim p_0 $, the number of kernels $K$, the learning rates $\rho_{\alpha}$ and $\rho_{\lambda}$, the number of iterations $L, J \in \mathbb{N}$, the total number of time steps $N_T$.
\vspace{0.3em}
\item[\textbf{while}] $n =0, 1, 2, \cdots, N_T-1$, \textbf{do} \\
\begin{description}
\item[-] \hspace{-0.25em} Propagate samples $\{x_i^{n}\}_{i=1}^N$ through the scheme \eqref{SDE:sample} to get $\{ \tilde{x}_i^{n+1}\}_{i=1}^{N}$. 
\item[-] \hspace{-0.3em}  Let $Y_n^{M_{t_n}} = p_n$ be the previous kernel learned filtering density. 
 Solve the time-inverse FBSDEs system \eqref{FBDSDEs:step} for $Y_{n+1}^{M_{t_n}}$ on spatial samples $\{ \tilde{x}_i^{n+1}\}_{i=1}^{N}$ through the iterative scheme \eqref{Scheme:Y:single}. $Y_{n+1}^{M_{t_n}}$ is the approximation for the predicted filtering density  $\tilde{p}(S_{t_{n+1}} \big| \mathcal{M}_{t_{n}})$. 
\item[-] \hspace{-0.25em}  Incorporate the observational information through Bayesian inference to get the updated filtering density $\tilde{p}(S_{t_{n+1}} \big| \mathcal{M}_{t_{n+1}})$ on the spatial samples $\{\tilde{x}_i^{n+1}\}_{i=1}^N$.   
\item[-] \hspace{-0.25em} Select kernel centers from the spatial samples $\{\tilde{x}_i^{n+1}\}_{i=1}^N$ by using the updated filtering density values on those samples. 
\item[-] \hspace{-0.75em}  Consider the approximated filtering density values  $\big\{\tilde{p}(S_{t_{n+1}} = \tilde{x}_i^{n+1} \big| \mathcal{M}_{t_{n+1}})\big\}_{i=1}^N$ as training data. Use the optimization procedure \eqref{Kernel-Learning} to obtain the kernel learned filtering density $p_{n+1}$ introduced in \eqref{KL-Pn}; 
\item[-] \hspace{-0.25em}  Carry out the resampling procedure to generate new samples $\{x_i^{n+1}\}_{i=1}^N$ that follow the kernel learned filtering density $p_{n+1}$;  
\end{description}
\item[\textbf{end while}]
\end{algorithmic}
\vspace{-1.2em}
\end{spacing}\\
\hline
\end{tabular}
\end{table}

One may notice that the state spatial samples that we use in the backward SDE filter have similar behavior to the particles in the particle filter method, which roughly characterize the filtering density in the nonlinear filtering problem. We want to emphasize that the backward SDE filter can also approximate the filtering density values on spatial samples. On the other hand, the particle filter only utilizes particle positions to construct an empirical distribution for the target state. Therefore, each spatial sample in the backward SDE filter carries more information about the state distribution than a particle in the particle filter. 
Moreover, since the kernel learned filtering density is a global continuous approximation for the state distribution, it covers wide range in the state space, which can provide more robust/stable performance for the backward SDE filter.


\section{Numerical experiments}\label{Numerics}
In this section, we use three numerical examples to demonstrate the performance of our kernel learning backward SDE filter method (BSDEF). In the first example, we focus on the BSDEF, and we use a synthetic nonlinear filtering problem to show mathematical behaviors of the BSDEF. In the second example, we solve the Lennard-Jones potential tracking problem, in which a target atom is moved by the intermolecular force generated by the Lennard-Jones potential. This is a benchmark problem in microphysics, and it has wide applications in material sciences. To demonstrate the effectiveness of the BSDEF, we compare the estimation performance of the BSDEF with the auxiliary particle filter (APF) \cite{APF} and the ensemble Kalman filter (EnKF) \cite{EvensenBook}, which are both state-of-the-art Bayes filter methods.
In the third example, we solve a Lorenz-96 tracking problem. The Lorenz-96 model forms the fundamental mathematical element for atmospheric data assimilation. It is well-known that the dynamical system driven by the Lorenz-96 dynamics is chaotic, and it's difficult to estimate the state of a Lorenz-96 model -- especially in high dimensional spaces. In this Lorenz-96  tracking example, we conduct systematic comparison studies, and we will show that the BSDEF is more accurate and more efficient compared with the APF and the EnKF.

The CPU that we use to implement all the numerical experiments is a 2.5 GHz Dual-Core Intel Core i7 processor with $16$ GB 2133 MHz LPDDR3 memory.

\subsection{A synthetic example}

In the first numerical example, we consider the following dynamical system:
\begin{equation}\label{Ex:synthetic}
\begin{aligned}
d S_t =& b(S_t) dt + {\bf \sigma} dW_t,\\
\end{aligned}
\end{equation}
where the target state $S: = ( X_1, X_2 )^{T}$ is a two dimensional vector driven by the dynamics $b$, which is defined as follows:
$$b(S) = \alpha \Big( \sin(X_2) +  \f{X_1}{1 + X_1}, \  \cos(X_1) +  \f{X_2}{1 + X_2} \Big)^T.$$
The observational data that we collect to estimate the state $S$ are direct observations, which are perturbed by Gaussian  noises with standard deviation $R$.

\begin{figure}[h!]
\begin{center}
\subfloat[Estimation performance of BSDEF]{\includegraphics[scale = 0.7]{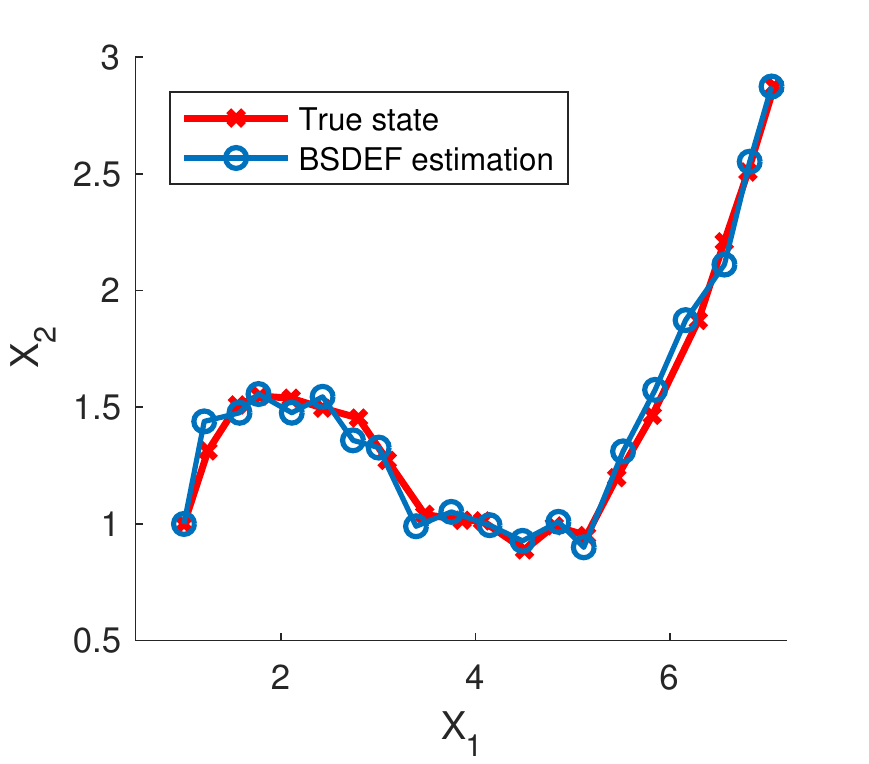} } 
\subfloat[Demonstration for kernel centers]{\includegraphics[scale = 0.7]{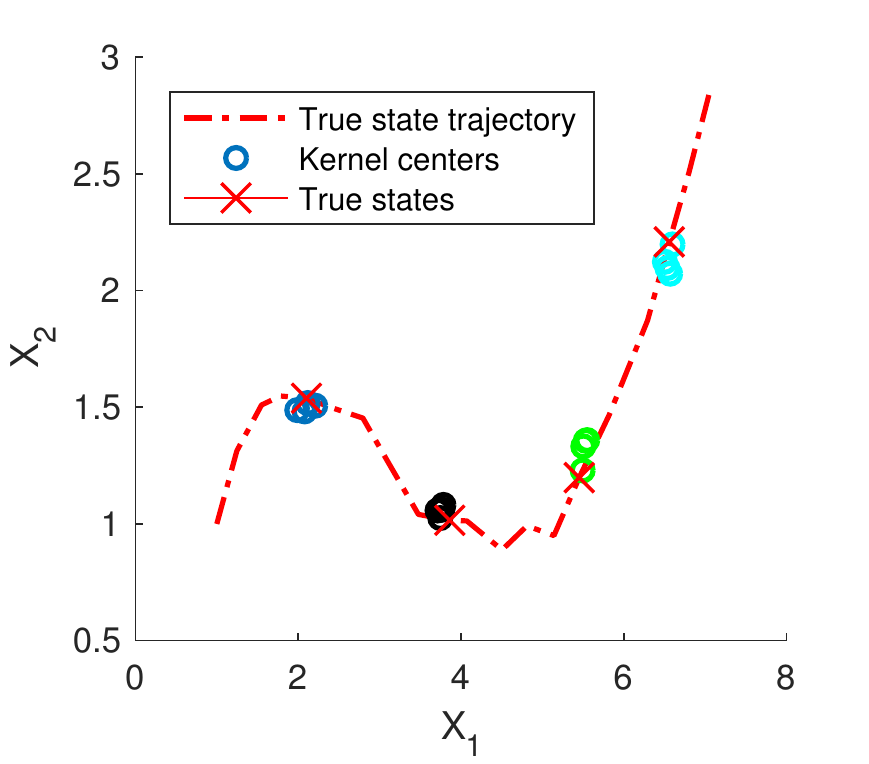} } \\
\subfloat[Confidence band for $X_1$ estimation]{\includegraphics[scale = 0.7]{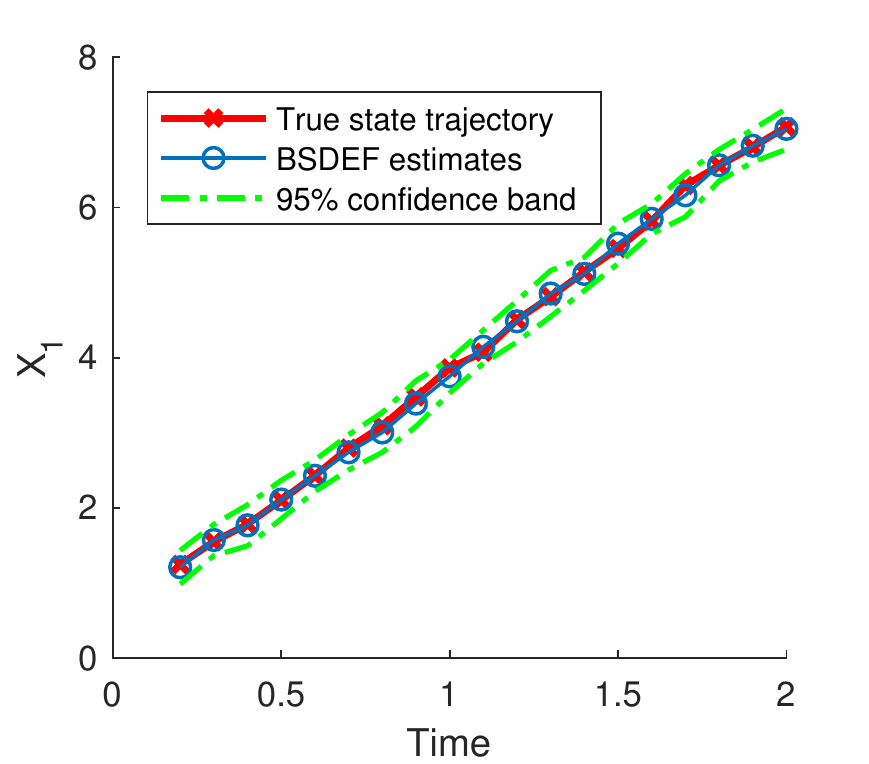} } 
\subfloat[Confidence band for $X_2$ estimation]{\includegraphics[scale = 0.7]{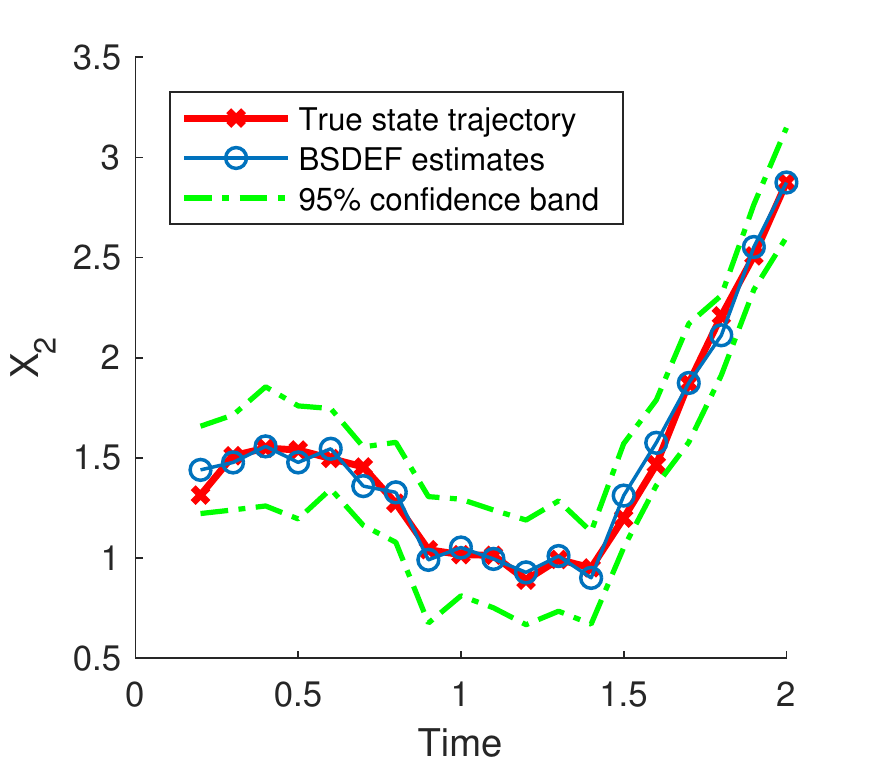} } 
\end{center}
\caption{Performance of BSDEF in solving the synthetic example}\label{Syn_Est} 
\end{figure}

In figure \ref{Syn_Est}, we present the performance of the BSDEF method in tracking the target state $S$ over the time period $[0, 2]$ with time step-size $\Delta t = 0.1$, i.e., $N_T = 20$, and we let $\alpha = 2$, $\sigma = 0.2 I_2$, and $R = 0.05 I_2$ in the nonlinear filtering problem.
The FBSDEs system is solved on $500$ spatial sample points, and we use $4$ Gaussian kernels to describe the filtering density. The fixed-point iterative scheme for solving the backward SDE is carried out with $10$ iteration steps, and the stochastic gradient descent optimization procedure for kernel learning is carried out with $100$ iteration steps. 
The initial guess for kernel weights is $0.5$ for each kernel, the initial guess for covariances is $2 I_{2}$ for each Gaussian kernel, and the learning rate is set to be $10^{-2}$. In Figure \ref{Syn_Est} (a), we present the state estimation performance. We can see from this subplot that our BSDEF gives very accurate estimates for the true target state trajectory. In Figure \ref{Syn_Est} (b), we plot the locations of Gaussian kernel centers at time steps $t = 0.5, 1, 1.5$, and $1.8$ by using blue, black, green and cyan circles, respectively, and the real target states at time steps $t=0.5, 1, 1.5, 1.8$ are given by red crosses. From this subplot, we can see that the kernel centers always surround the true state of the target, which guarantees that the high probability regions of the filtering density cover the target state. To show more detailed behavior of the filtering density obtained by the BSDEF, we plot $95\%$ marginal confidence bands of the filtering density for $X_1$ and $X_2$ in subplots Figure \ref{Syn_Est} (c) and Figure \ref{Syn_Est} (d), respectively, where the confidence bands are plotted by green dashed curves. From these subplots, we can see that the true target state is always within the $95\%$ confidence bands of the filtering density obtained by our BSDEF method. 

An important concept in the backward SDE approach for nonlinear filtering problems is to use backward SDEs to generate filtering density values on scattered spatial sample points, and then use the locations of those samples together with their filtering density values to describe the conditional PDF of the target state. Therefore, the number of spatial sample points has strong influence to the accuracy of the BSDEF method. To demonstrate the convergence performance of BSDEF with respect to the number of spatial samples, we solve the nonlinear filtering problem \eqref{Ex:synthetic} repeatedly $100$ times by using different random seeds to generate samples of random variables. 
\begin{figure}[h!]
\begin{center}
\includegraphics[scale = 0.8]{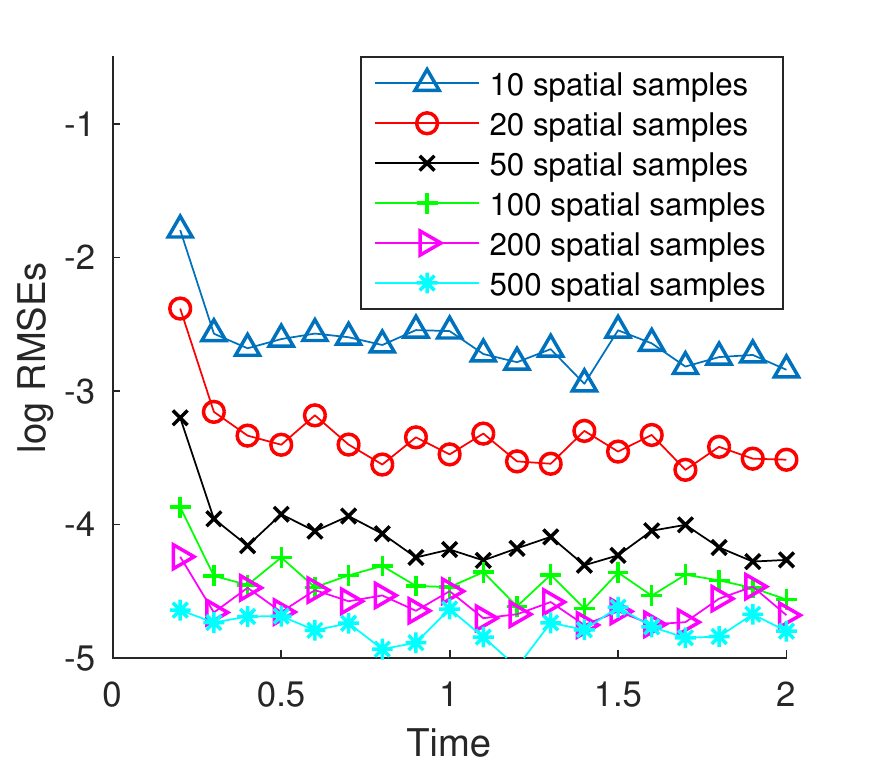} 
\end{center}
\caption{Convergence of BSDEF estimates with respect to spatial samples.}\label{Syn_Convergence} 
\end{figure}
In Figure \ref{Syn_Convergence}, we plot the root mean square errors (RMSEs) between the true state and our estimated state with respect to time by using $10$, $20$, $50$, $100$, $200$, and $500$ spatial samples in the spatial approximation for the BSDEF, where the horizontal axis is the time axis and the vertical axis shows the natural logarithms of RMSEs. From this figure, we can see clearly that the BSDEF gives more and more accurate estimates for the target state as the number of spatial samples increases. 

\subsection{The Lennard-Jones potential tracking problem}

In this example, we track a target atom driven by the Lennard-Jones potential $V_{LJ}$, which is an important intermolecular potential model that has been most extensively studied and applied.  It is considered as the archetype model for intermolecular interactions. The task of tracking the trajectory of an atom driven by the Lennard-Jones potential through direct observations received by scanning transmission electron microscopy is also one of the key mathematical challenges in a recently developed novel technique called the ``atomic forge'', which aims to design and assemble materials atom-by-atom \cite{Kalinin-Atom, Bao_Atomic_2021}. 

The mathematical model of the Lennard-Jones potential in its AB form is given by $\ds V_{LJ} = \f{A}{r^{12}} - \f{B}{r^6}$, where $A$, $B$ are two energy parameters that determine the physical properties of the intermolecular potential and $r$ is the distance between two atoms. In this work, we track a freely moving target atom and fix the location of the other one, which is called the platform atom in this example. 
The force between these two atoms is the gradient of $V_{LJ}$, i.e. $\nabla V_{LJ}$. 
When these two atoms are at moderate distance from each other, the intermolecular force $\nabla V_{LJ}$ appears attractive and it drags the target atom towards the fixed platform atom. On the other hand, when these two atoms get too close, the atomic force $\nabla V_{LJ}$ becomes repulsive quickly, which pushes the moving target atom away from the fixed platform atom. The state equation that describes the moving target atom is given by
\begin{equation}\label{Ex:LJ}
dS_t = - \nabla V_{LJ}(S_t) dt + \sigma dW_t,
\end{equation}
where $S : = (x, y)^T$ is the 2D location of the target atom on a material surface, the stochastic noise term $\sigma dW_t$ models random movements of the target atom. When the moving target atom is close to the fixed platform atom, a random movement towards the fixed atom may result a very large repulsive force, which pushes the target atom away. Since such an intermolecular force is unpredictable, and it could dramatically change the location of the target atom, the Lennard-Jones potential problem is a challenging benchmark experiment to examine performances of different nonlinear filtering methods.

In our numerical experiment, we track the location of the target atom over the time period $[0, 30]$ with observational gap $\Delta t = 0.3$, i.e. $N_T = 100$. The observational data we collect are target positions, which are perturbed by Gaussian noises with standard deviation $R = 0.01 I_2$.  We let the fixed platform atom be at the origin of the $xy$-plan, $A = 16$, $B=4$ and $\sigma = 0.02 I_2$ in the dynamical model \eqref{Ex:LJ}, and the initial position of the target atom is given at  $(1.8, 2.2)^T$.

\begin{figure}[h!]
\begin{center}
\subfloat[Comparison of estimation for $x$]{\includegraphics[scale = 0.6]{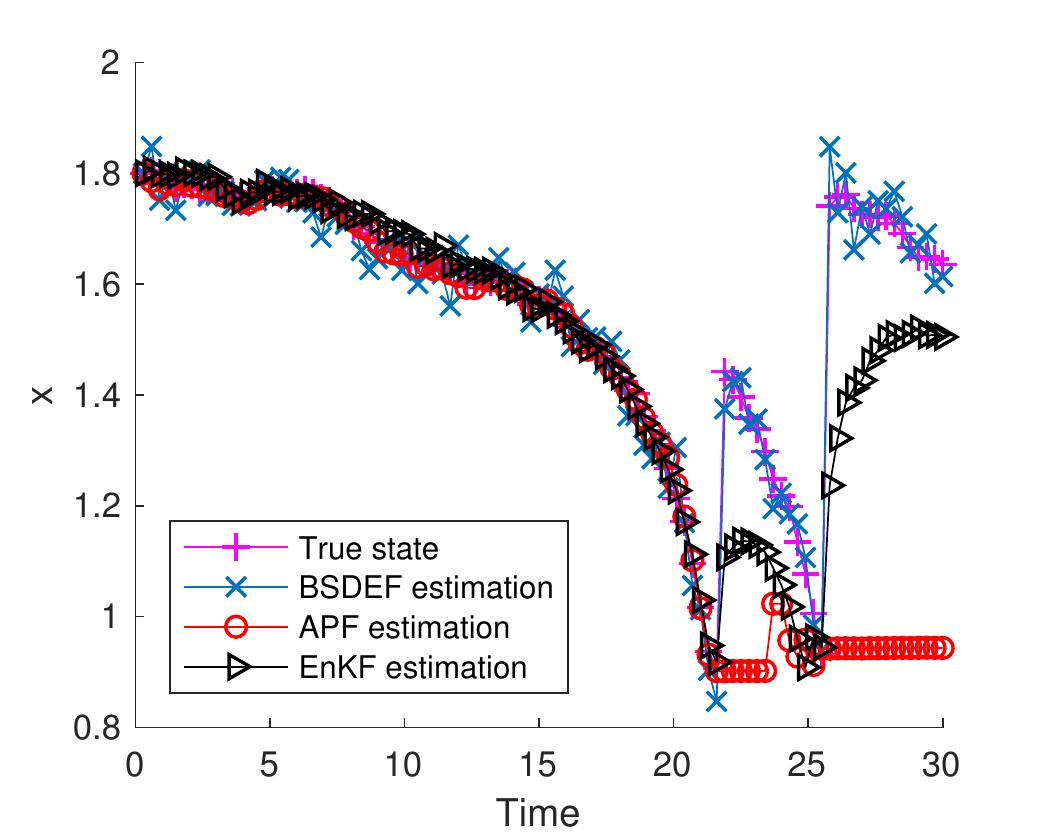} } 
\subfloat[Comparison of estimation for $y$]{\includegraphics[scale = 0.6]{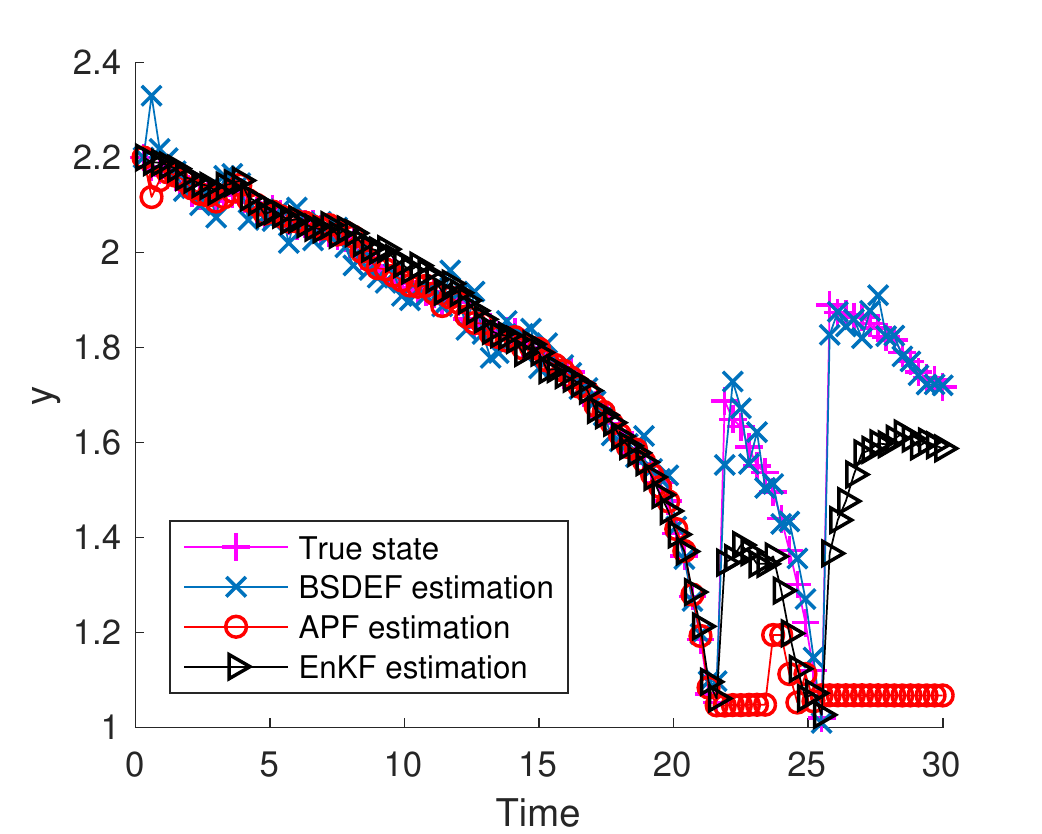} }
\end{center}
\caption{Comparison of estimation performances}\label{LJ_comparison} 
\end{figure}
In Figure \ref{LJ_comparison}, we present the comparison of estimation performances between our BSDEF, the auxiliary particle filter (APF) and the Ensemble Kalman filter (EnKF), where the subplot (a) shows the estimation performances in tracking the $x$-component of $S$, and the subplot (b) shows the estimation performances in tracking the $y$-component of $S$. To carry out the BSDEF, we use $200$ spatial sample points and $4$ Gaussian kernels to learn the filtering density. For the APF, we use $500$ particles to describe the filtering density with $10$ auxiliary Monte Carlo samples, and the EnKF is implemented by using an ensemble of $1000$ realizations of Kalman filter samples. The true target atom (the magenta curve marked by pluses) moves towards the fixed atom at the origin, and it was pushed away by the intermolecular force twice after the time instants $t = 20$ and $t=25$. Throughout of the tracking period, the BSDEF estimation (the blue curve marked by crosses) always accurately captures the true locations of the target atom. On the other hand, the APF (the red curve marked by circles) fails after the first (unpredicted) swift location change due to highly nonlinear behaviors of the target atom and the degeneracy of the particles, and the EnKF (the black curve marked by triangles) loses the target after swift location changes as well. From this experiment, we can see that our BSDEF method has more robust performance in atom tracking compared with the APF and the EnKF.

\subsection{The Lorenz-96 tracking problem}
In the third numerical example, we carry out comprehensive comparisons between the BSDEF, the APF and the EnKF. The state equation that we consider in this example is driven by the Lorenz-96 dynamics, i.e.,
\begin{equation}\label{Ex:Lorenz}
dS_t = {\bf b}(S_t) dt + \sigma dW_t,
\end{equation}
where $S = (x_{1}, x_2, \cdots, x_d)$ is a $d$-dimensional state vector. The Lorenz-96 dynamics ${\bf b}(S) = (b_1(S), b_2(S), \cdots, b_d(S))^T$ is defined by
\begin{equation*}
b_{i}(S) = (x_{i+1} - x_{i-2})x_{i-1} - x_{i} + F, \quad i = 1, 2, \cdots, d,
\end{equation*}
where $x_{-1} = x_{d-1}$, $x_{0} = x_{d}$ and $x_{d+1} = x_{1}$. It is well-known that when the forcing constant $F=8$, the state dynamics \eqref{Ex:Lorenz} has chaotic behaviors, hence tracking the state of Lorenz-96 model is a challenging benchmark nonlinear filtering problem. In addition, the dimension of the Lorenz-96 dynamics can be arbitrarily chosen ($d\geq 4$), which makes the Lorenz-96 tracking problem an ideal example to test the high dimensional estimation performance of nonlinear filtering methods.

In this Lorenz-96 tracking problem, we use the BSDEF, the APF, and the EnKF to estimate the state of $S$ over the time period $[0, 1]$ with time step $\Delta t = 0.02$, i.e. $N_T = 50$, and we let $\sigma = 0.1 I_d$. The initial true state $S_0$ is chosen as $2 + 4 \gamma$, where $\gamma$ is a standard $d$-dimensional Gaussian random variable. The initial guess for the state is $S_0$ with a noisy perturbation, and we let the noise be generated by a Gaussian random variable with standard deviation $0.5 I_d$.

In the first numerical experiment, we let the dimension of the Lorenz-96 model be $d = 10$ and assume that direct observations are available, i.e. $M_{t_n} = S_{t_n} +  \xi$, where $\xi$ is a standard Gaussian random variable with standard deviation $0.1 I_d$, which represents the observational error, and $\xi$ is independent of the Brownian motion $W$ in \eqref{Ex:Lorenz}.  To carry out this comparison experiment, we select $800$ spatial samples in the BSDEF, and we use $10$ Gaussian kernels to learn the filtering density. The fixed-point iterative scheme for solving the backward SDE is carried out with $10$ iteration steps, and the stochastic gradient descent optimization procedure for kernel learning is carried out with $100$ iteration steps. 
The initial guess for the kernel weight is $1$ for each kernel, the initial guess for the covariance is $I_{d}$ for each kernel, and the learning rate is set to be $10^{-2}$. 
For the APF, we use $2000$ particles to describe the distribution of the target state with $10$ auxiliary Monte Carlo samples, and we use $3000$ realizations of Kalman filter samples in the EnKF. 
\begin{figure}[h!]
 \vspace{-0.5em}
\begin{center}
\subfloat[Estimates for $x_3$]{\includegraphics[scale = 0.49]{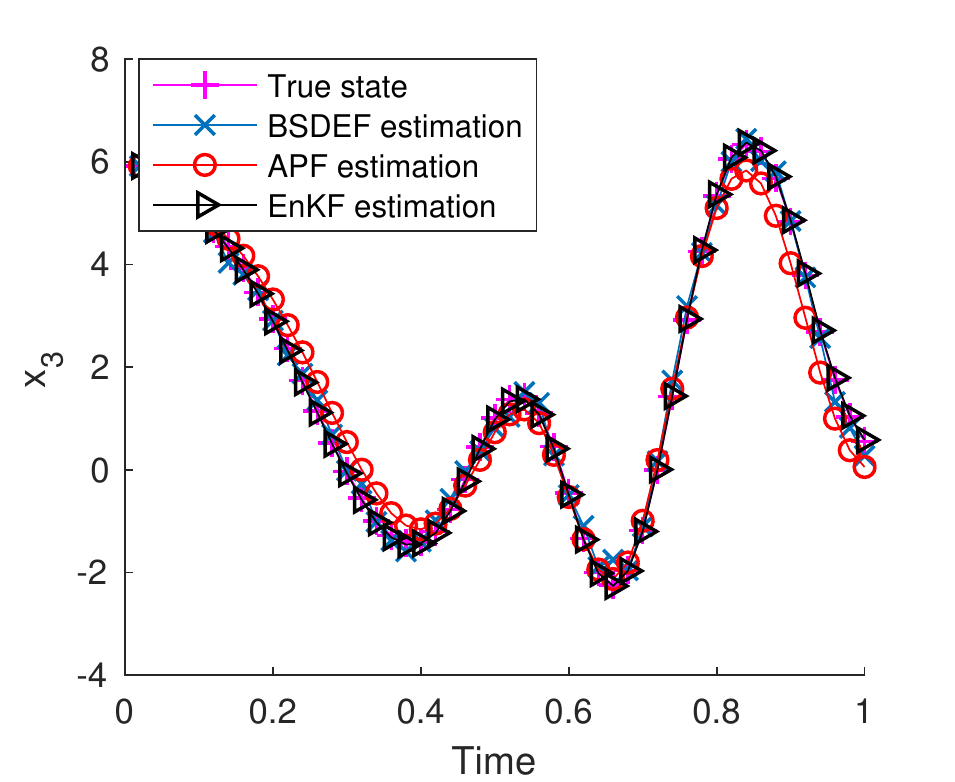} }  \hspace{-0.5em}
\subfloat[Estimates for $x_6$]{\includegraphics[scale = 0.49]{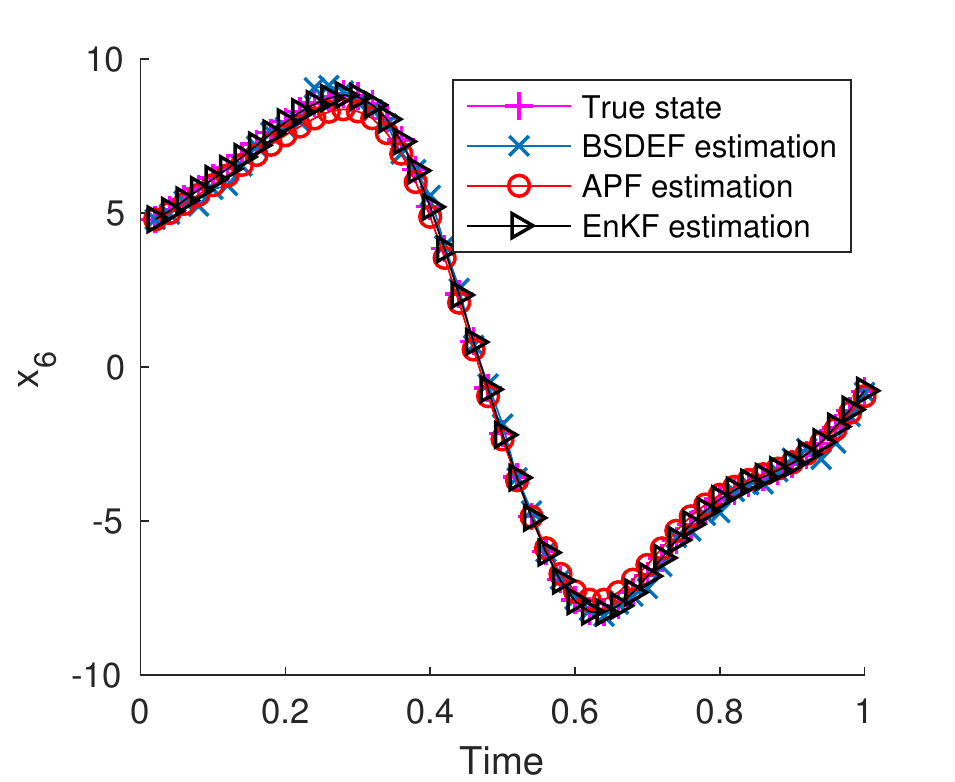} } \hspace{-0.5em}
\subfloat[Estimates for $x_{9}$]{\includegraphics[scale = 0.49]{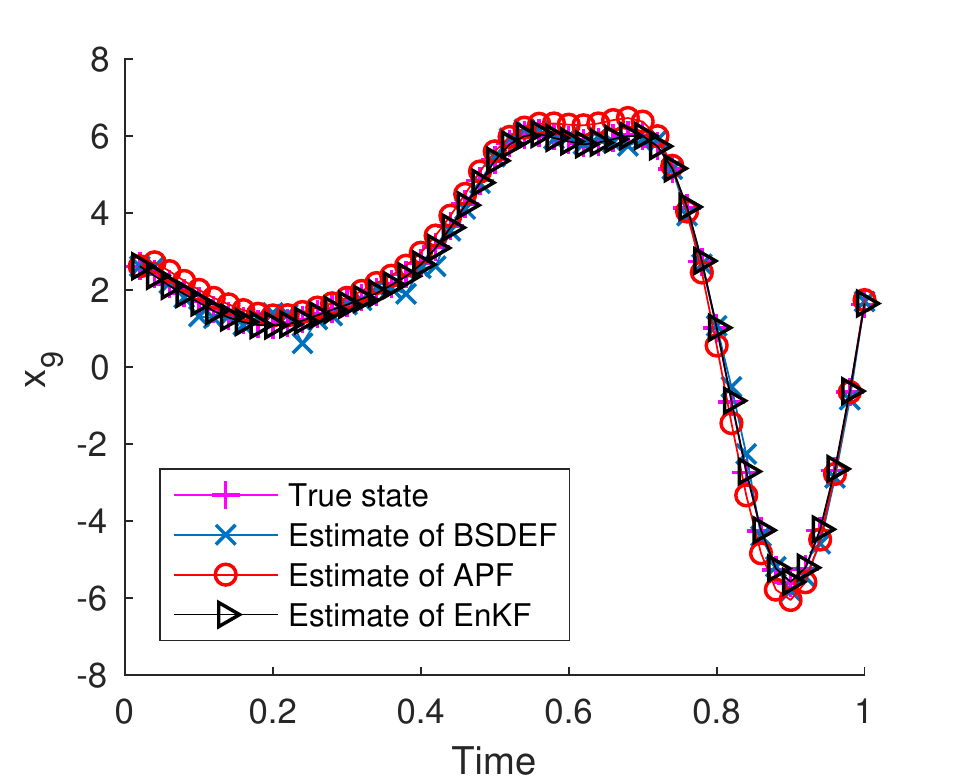} }
\end{center}
\caption{Estimation performance with linear observations}\label{Lorenz_Linear} 
\end{figure}
In Figure \ref{Lorenz_Linear}, we present performances of the BSDEF, the APF and the EnKF in estimating the target state in spatial dimensions $x_3$, $x_6$ and $x_9$. We can see from this figure that all three methods work very well in this experiment.

However, in the following experiment, we replace the direct linear observation function by the cubic root observation function, i.e. $M_{t_n} = (S_{t_n})^{\f{1}{3}} + \xi$, and we keep all the other parameters unchanged (both for the nonlinear filtering problem and for the filtering methods) as in the previous experiment.
\begin{figure}[h!]
 \vspace{-0.5em}
\begin{center}
\subfloat[Estimates for $x_3$]{\includegraphics[scale = 0.49]{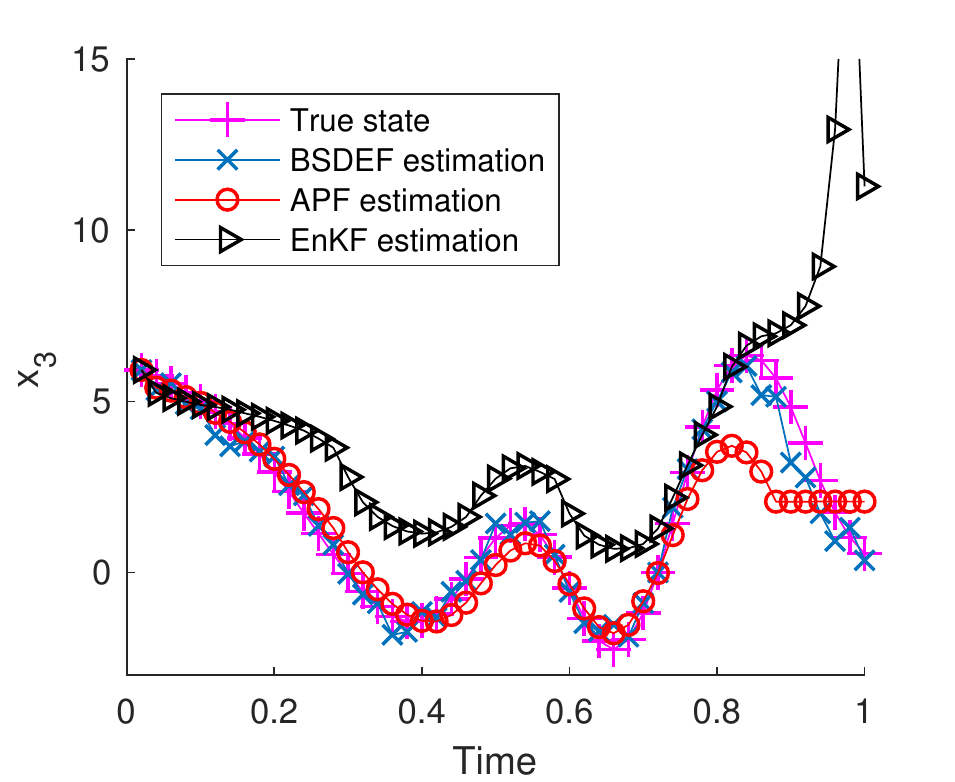} }  \hspace{-.5em}
\subfloat[Estimates for $x_6$]{\includegraphics[scale = 0.49]{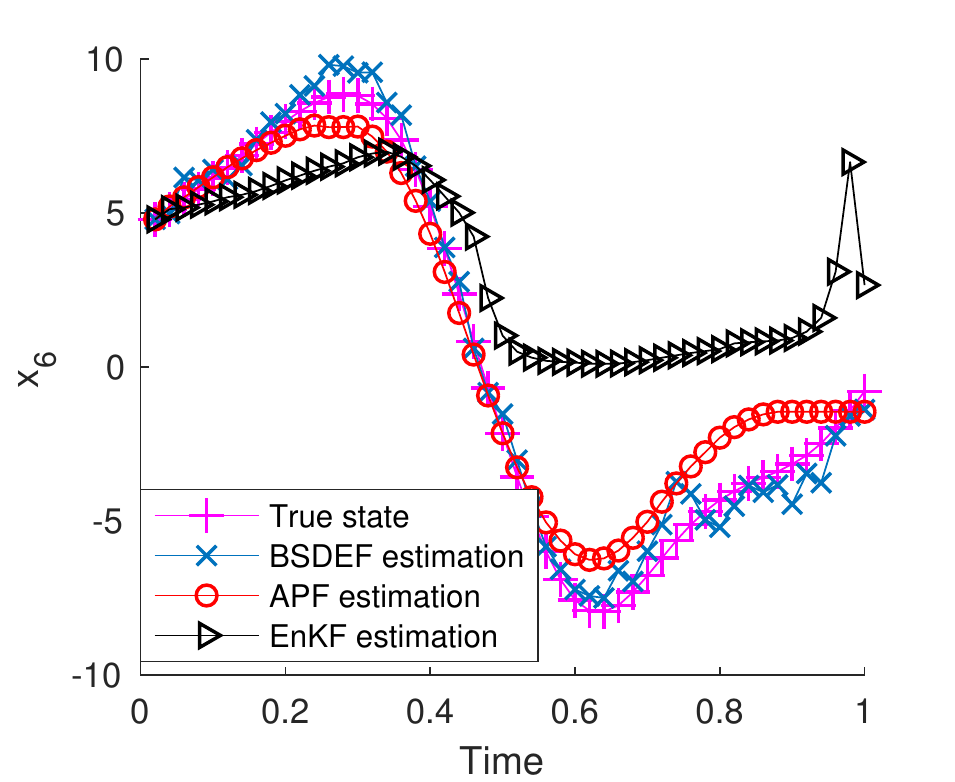} } \hspace{-.5em}
\subfloat[Estimates for $x_{9}$]{\includegraphics[scale = 0.49]{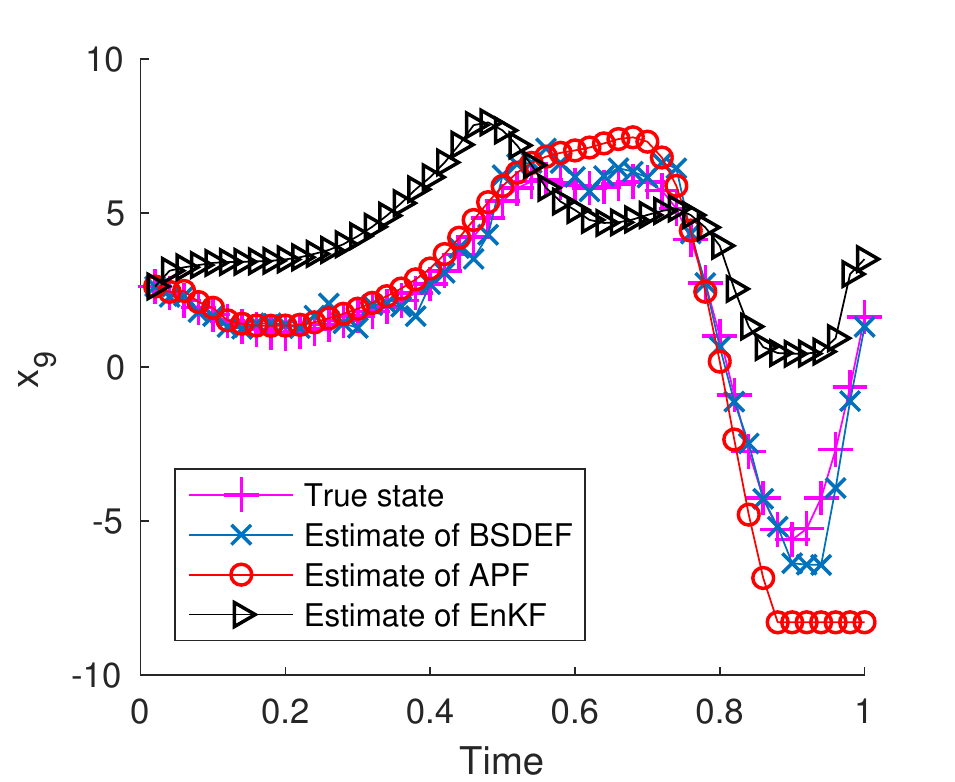} }
\end{center}
\caption{Estimation performance with cubic observations}\label{Lorenz_Cubic} 
\end{figure}
In Figure \ref{Lorenz_Cubic}, we present performances of all the three methods in estimating $x_3$, $x_6$, and $x_9$, where the true state trajectories in three subplots are magenta curves marked by pluses, the BSDEF estimates are blue curves marked by crosses, the APF estimates are red curves marked by circles, and the EnKF estimates are black curves marked by triangles. With nonlinear observations in this experiment, although the true state has exactly the same trajectory as presented in Figure \ref{Lorenz_Linear},  the EnKF can not provide accurate estimates for the target state. The APF gives  accurate estimates at beginning. However, as we tracking the target for more steps, the particle filter starts to suffer from the degeneracy problem, which may also be caused by indirect observations, and the estimation performance deteriorates quickly.
On the other hand, our BSDEF method is always on track and it accurately captures the true state of the target. 

To further demonstrate the estimation performance, we repeat the above nonlinear observation experiment $50$ times by using different seeds to generate random numbers for the random variables in the nonlinear filtering problem and calculate the RMSEs with respect to tracking time. 
\begin{figure}[h!]
 \vspace{-0.5em}
\begin{center}
\includegraphics[scale = 0.65]{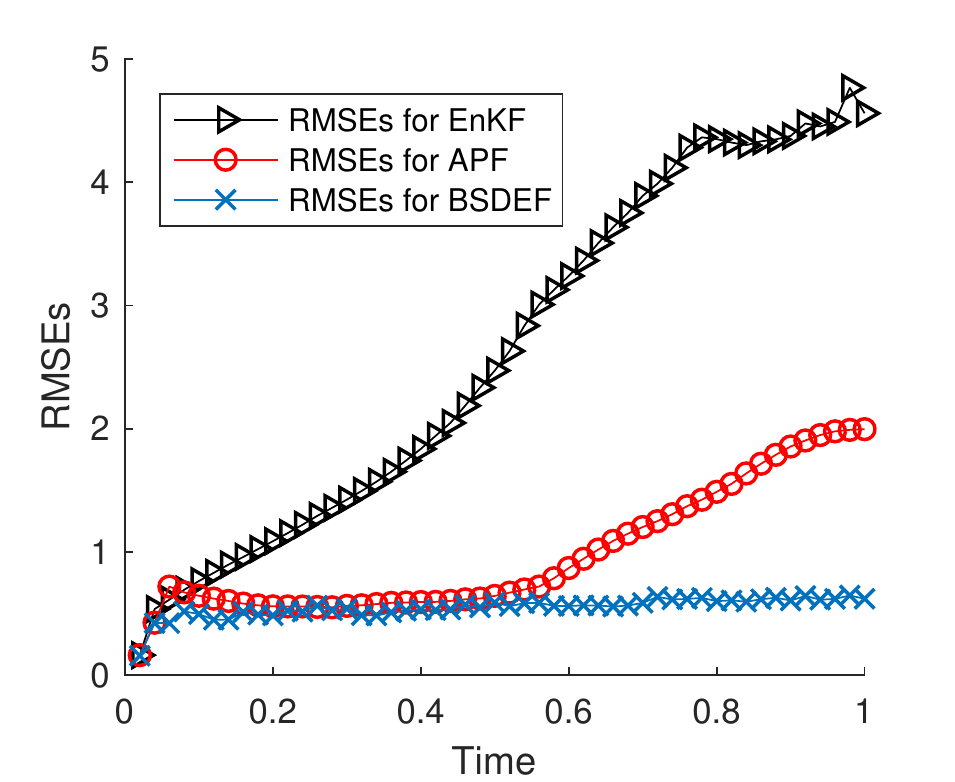} 
\end{center}
\caption{Comparison of RMSEs for $d = 10$}\label{Lorenz_10d_RMSE} 
\end{figure}
The accuracy of three methods is presented in Figure \ref{Lorenz_10d_RMSE}, where the black curve marked by triangles shows the RMSEs of the EnKF, the red curve marked by circles shows the RMSEs of the APF, and the blue curve marked by crosses is the RMSEs of the BSDEF.
We can see from this figure that the BSDEF has the lowest errors among all three methods, and it has the most stable performance with respect to the tracking time. The APF has comparable performance with the BSDEF at beginning. But the RMSEs of the APF start to grow near the time instant $t=0.6$. Therefore, the long term performance of the APF is not as good as the BSDEF. On the other hand, the EnKF has the highest RMSEs. The average CPU time to implement the BSDEF in this experiment is $5.71$ seconds. The average CPU time to implement the APF is $14.43$ seconds, and the average CPU time to implement the EnKF is $6.67$ seconds. For this $10$ dimensional experiment, we can see that the BSDEF outperforms the APF and the EnKF in both accuracy and efficiency.

In order to make the advantageous performance of our kernel learning backward SDE filter more convincing, we carry out the above RMSEs comparison experiment on higher dimensional Lorenz-96 tracking problems, i.e. $d=15$ and $d=20$. For the $d=15$ case, we use $1000$ spatial samples in the BSDEF with $15$ Gaussian kernels, $4000$ particles in the APF, and $5000$ realizations of Kalman filter samples in the EnKF.  For the $d=20$ case, we use $1500$ spatial samples for the BSDEF with $20$ Gaussian kernels, $6000$ particles in the APF, and $10,000$  realizations of Kalman filter samples in the EnKF.   
\begin{figure}[h!]
 \vspace{-0.5em}
\begin{center}
\subfloat[Comparison of RMSEs for $d = 15$]{\includegraphics[scale = 0.65]{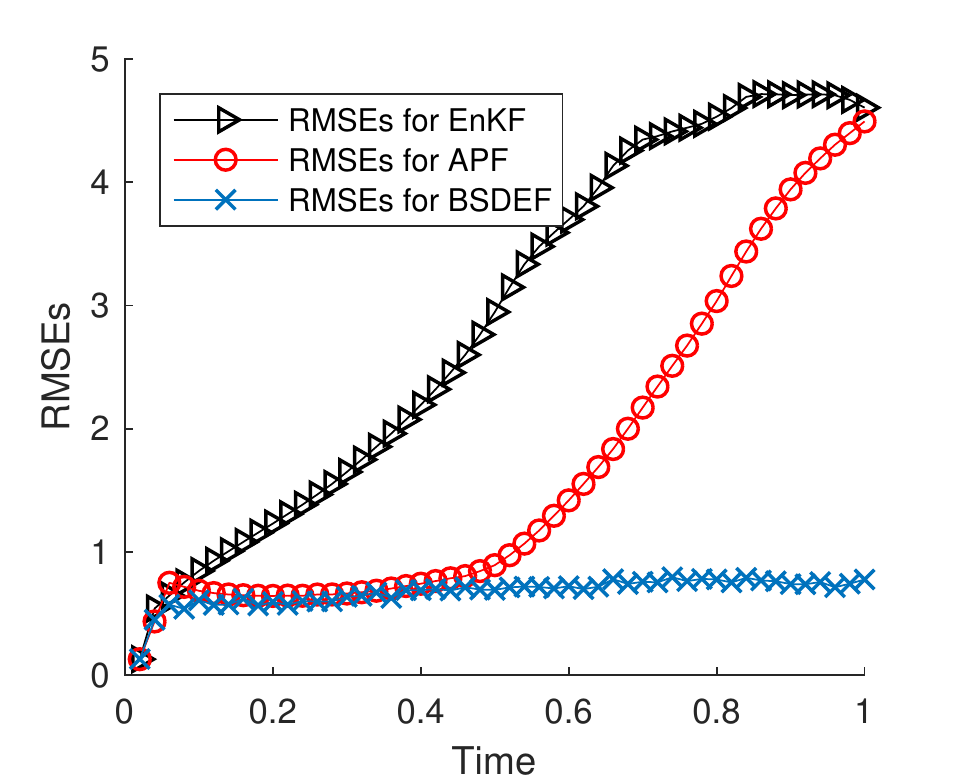} }
\subfloat[Comparison of RMSEs for $d = 20$]{\includegraphics[scale = 0.65]{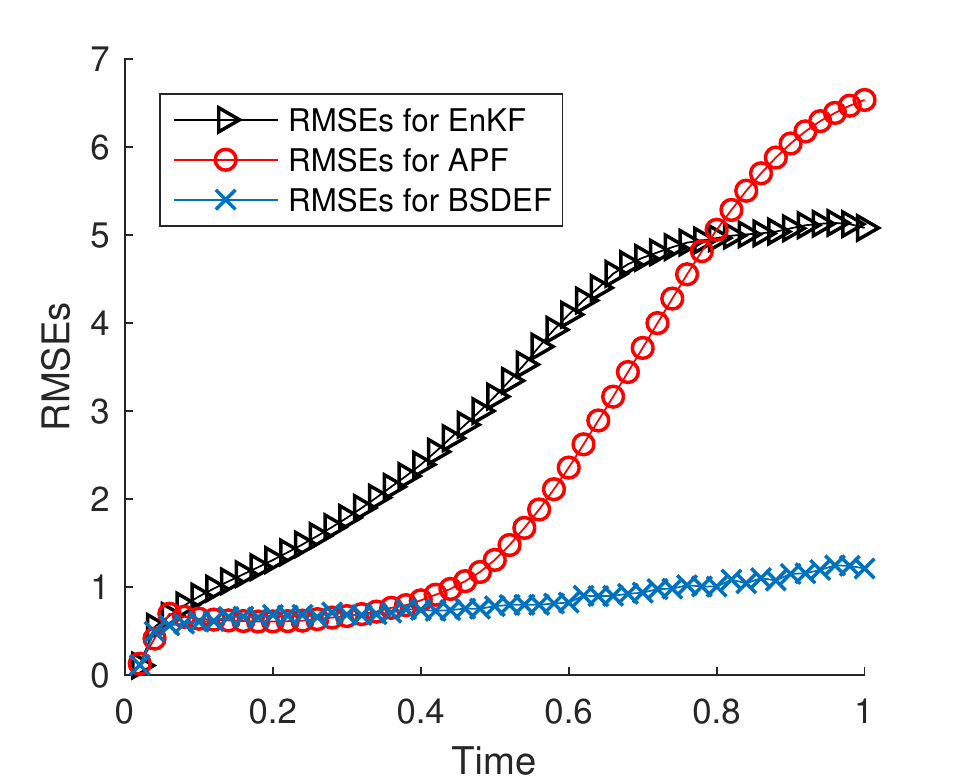} }
\end{center}
\caption{Comparison of RMSEs for $d=15$ and $d = 20$}\label{RMSEs} 
\vspace{-0.5em}
\end{figure}
The RMSEs are presented in Figure \ref{RMSEs}. We can see from this figure, as well as the $10$-dimensional comparison result, that the BSDEF has the lowest and the most stable RMSEs. The APF always has comparable performances with the BSDEF at beginning, and its RMSEs grows as more and more estimation steps are carried out. In all the three RMSE experiments, the EnKF has high errors. However, we notice that both the BSDEF and the EnKF maintain their estimation accuracy when the dimension of the problem increases. On the other hand, the RMSEs of the APF increase as the dimension of the problem increases. In the $20$-dimensional experiment, the APF starts to produce higher RMSEs compared with EnKF from the time instant $t = 0.8$. This shows that the particle filter has poor performance for high dimensional problems and for long term estimations. 

In Table \ref{table}, we summarize the RMSE experiments for $d=10$, $d=15$ and $d = 20$. 
\begin{table}[h!]\footnotesize
\caption{Summary of numerical comparison in the Lorenz-96 model tracking.} \label{table}
\vspace{-1em}
\begin{center}
\begin{tabular}{|c|c|c|c|c|c|c|c|c|c|c|}
\hline
\multirow{2}{*}{} & \multicolumn{3}{c|}{d = 10} & \multicolumn{3}{c|}{d = 15} & \multicolumn{3}{c|}{d = 20} \\ \cline{2-10} 
     &      Samples   &    Time      &  RMSEs &  Samples    &     Time      &   RMSEs  &  Samples    &    Time      &  RMSEs    \\
\hline    {\bf BSDEF }  &     800    &   {\bf  5.71 }   &  {\bf 3.92} & 1,000 &{\bf 9.74}  &{\bf 4.77 } & 1,500 &{\bf 18.57}  & {\bf 5.40 }  \\
\hline     APF   &       2,000    &     14.43     &  7.79 & 3,000 & 24.12 &15.17  & 6,000 &56.34  & 23.51   \\
\hline     EnKF   &      3,000    &   6.67  &  21.21 & 5,000 & 12.31 & 23.13  & 10,000 & 27.14  & 25.20   \\
\hline  \multicolumn{10}{c}{$\ast$ The unit for CPU Time is second. RMSEs are accumulated RMSEs over the tracking period. } 
\end{tabular}\end{center}\vspace{-2em}
\end{table}
The results presented in this table again verifies that the BSDEF is most accurate and most efficient among all three state-of-the-art methods.

\section{Conclusions}\label{Conclusion}

In this paper, we introduced a kernel learning backward SDE filter to solve the nonlinear filtering problem. The main theme of our kernel learning approach is to treat the discrete filtering density values obtained by the backward SDE filter as simulation data, and then we use kernel learning to learn a continuous global approximation for the entire filtering density from the simulation data. The primary advantage of the kernel learned filtering density is that it provides a comprehensive description for the filtering density in the entire state space, which makes the kernel learning backward SDE filter more accurate and more stable in estimating the target state. Numerical experiments are presented to demonstrate the effectiveness and efficiency of kernel learning backward SDE filter in solving a synthetic problem and two benchmark application problems.

\vspace{1em}

\subsection*{References.}

\bibliographystyle{plain}


\begin{thebibliography}{10}

\bibitem{MCMC-PF}
C.~Andrieu, A.~Doucet, and R.~Holenstein.
\newblock Particle markov chain monte carlo methods.
\newblock {\em J. R. Statist. Soc. B}, 72(3):269--342, 2010.

\bibitem{BaoCC_2019}
F.~Bao, Y.~Cao, and H.~Chi.
\newblock Adjoint Forward Backward Stochastic Differential Equations Driven
by Jump Processes and Its Application to Nonlinear Filtering Problems.
\newblock {\em Int. J. Uncertain Quantif}, 9(2):143-159, 2019.

\bibitem{BaoCH_CiCP}
F.~Bao, Y.~Cao, and X.~Han.
\newblock An Implicit Algorithm of Solving Nonlinear Filtering Problems.
\newblock {\em Commun Comput Phys}, 16(2):382-402, 2014.

\bibitem{Bao_CMS}
F.~Bao, Y.~Cao, and X.~Han.
\newblock Forward backward doubly stochastic differential equations and optimal
  filtering of diffusion processes.
\newblock {\em Commun Math Sci .}, 18(3):635--661, 2020.

\bibitem{Bao_first}
F.~Bao, Y.~Cao, A.~J. Meir, and W.~Zhao.
\newblock A first order scheme for backward doubly stochastic differential
  equations.
\newblock {\em SIAM/ASA J. Uncertain. Quantif.}, 4(1):413--445, 2016.

\bibitem{Bao_Zakai}
F.~Bao, Y.~Cao, C.~Webster, and G.~Zhang.
\newblock A hybrid sparse-grid approach for nonlinear filtering problems based
  on adaptive-domain of the {Z}akai equation approximations.
\newblock {\em SIAM/ASA J. Uncertain. Quantif.}, 2(1):784--804, 2014.

\bibitem{BCZ_2011}
F.~Bao, Y.~Cao, and W.~Zhao.
\newblock Numerical Solutions for Forward Backward Doubly Stochastic
Differential Equations and Zakai Equations.
\newblock {\em Int. J. Uncertain Quantif}, 4(1):351-367, 2011.


\bibitem{BCZ_2015}
F.~Bao, Y.~Cao, and W.~Zhao.
\newblock A First Order Semi-discrete Algorithm for Backward Doubly
Stochastic Differential Equations
\newblock {\em Discrete Contin. Dyn. Syst. Ser. B}, 5(2):1297-1313, 2015.

\bibitem{BCZ_2018}
F.~Bao, Y.~Cao, and W.~Zhao.
\newblock A backward doubly stochastic differential equation approach for
  nonlinear filtering problems.
\newblock {\em Commun. Comput. Phys.}, 23(5):1573--1601, 2018.


\bibitem{BSDE_filter}
F.~Bao and V.~Maroulas.
\newblock Adaptive meshfree backward {SDE} filter.
\newblock {\em SIAM J. Sci. Comput.}, 39(6):A2664--A2683, 2017.

\bibitem{SGD_NIPS2007}
L. Bottou and O. Bousquet.
\newblock The tradeoffs of large scale learning.
\newblock In J.~Platt, D.~Koller, Y.~Singer, and S.~Roweis, editors, {\em
  Advances in Neural Information Processing Systems}, volume~20. Curran
  Associates, Inc., 2008.

\bibitem{CT1}
A.~J. Chorin and X.~Tu.
\newblock Implicit sampling for particle filters.
\newblock {\em Proc. Nat. Acad. Sc. USA}, 106:17249--17254, 2009.

\bibitem{Bao_Atomic_2021}
O.~Dyck, M.~Ziatdinov, S.~Jesse, F.~Bao, A.~Yousefzadi Nobakht, A.~Maksov, B.G.
  Sumpter, R.~Archibald, K.J.H. Law, and S.V. Kalinin.
\newblock Probing potential energy landscapes via electron-beam-induced single
  atom dynamics.
\newblock {\em Acta Materialia}, 203:116508, 2021.

\bibitem{BSDE_finance}
N.~El~Karoui, S.~Peng, and M.~C. Quenez.
\newblock Backward stochastic differential equations in finance.
\newblock {\em Math. Finance}, 7(1):1--71, 1997.

\bibitem{EvensenBook}
G.~Evensen.
\newblock {\em Data assimilation: the ensemble {K}alman filter}.
\newblock Springer, 2006.

\bibitem{Gobet-Zakai}
E. Gobet, G. Pag{\`e}s, H. Pham, and J. Printems.
\newblock Discretization and simulation of the {Z}akai equation.
\newblock {\em SIAM J. Numer. Anal.}, 44(6):2505--2538 (electronic), 2006.

\bibitem{Multi-Kernel_Learning}
M. G\"{o}nen and E. Alpayd\i~n.
\newblock Multiple kernel learning algorithms.
\newblock {\em J. Mach. Learn. Res.}, 12:2211--2268, 2011.

\bibitem{particle-filter}
N.J Gordon, D.J Salmond, and A.F.M. Smith.
\newblock Novel approach to nonlinear/non-gaussian bayesian state estimation.
\newblock {\em IEE PROCEEDING-F}, 140(2):107--113, 1993.

\bibitem{Kernel_learning}
T. Hofmann, B. Sch\"{o}lkopf, and A.~J. Smola.
\newblock Kernel methods in machine learning.
\newblock {\em Ann. Statist.}, 36(3):1171--1220, 2008.

\bibitem{HU-Zakai}
Y.~Hu, G.~Kallianpur, and J.~Xiong.
\newblock An approximation for the {Z}akai equation.
\newblock {\em Appl. Math. Optim.}, 45(1):23--44, 2002.

\bibitem{UnKF}
S.J. Julier and J.K. Uhlmann.
\newblock Unscented filtering and nonlinear estimation.
\newblock {\em Proceedings of the IEEE}, 92:401--422, 2004.

\bibitem{Kalinin-Atom}
S.~Kalinin, A.~Borisevich, and S.~Jesse.
\newblock Fire up the atom forge.
\newblock {\em Nature}, 22 November 2016.

\bibitem{Kalman1961}
R.~E. Kalman and R.~S. Bucy.
\newblock New results in linear filtering and prediction theory.
\newblock {\em Transactions of the ASME--Journal of Basic Engineering},
  83(Series D):95--108, 1961.

\bibitem{Kang-PF}
K.~Kang, V.~Maroulas, I.~Schizas, and F.~Bao.
\newblock Improved distributed particle filters for tracking in a wireless
  sensor network.
\newblock {\em Comput. Statist. Data Anal.}, 117:90--108, 2018.

\bibitem{Kloeden1992}
P.~E. Kloeden and E. Platen.
\newblock {\em Numerical solution of stochastic differential equations},
  volume~23 of {\em Applications of Mathematics (New York)}.
\newblock Springer-Verlag, Berlin, 1992.

\bibitem{MTAC2012}
M.~Morzfeld, X.~Tu, E.~Atkins, and A.~J. Chorin.
\newblock A random map implementation of implicit filters.
\newblock {\em J. Comput. Phys.}, 231(4):2049--2066, 2012.

\bibitem{Pardoux1991}
{\'E}.~Pardoux and S.~Peng.
\newblock Backward stochastic differential equations and quasilinear parabolic
  partial differential equations.
\newblock In {\em Stochastic partial differential equations and their
  applications ({C}harlotte, {NC}, 1991)}, volume 176 of {\em Lecture Notes in
  Control and Inform. Sci.}, pages 200--217. Springer, Berlin, 1992.

\bibitem{PP1994}
{\'E}. Pardoux and S. Ge Peng.
\newblock Backward doubly stochastic differential equations and systems of
  quasilinear {SPDE}s.
\newblock {\em Probab. Theory Related Fields}, 98(2):209--227, 1994.

\bibitem{Peng_ICM}
S. Peng.
\newblock Backward stochastic differential equation, nonlinear expectation and
  their applications.
\newblock Proceedings of the ICM 2010, pages 393 -- 432, 2011.

\bibitem{APF}
M.~K. Pitt and N. Shephard.
\newblock Filtering via simulation: auxiliary particle filters.
\newblock {\em J. Amer. Statist. Assoc.}, 446(94):590--599, 1999.

\bibitem{Stochastic_Approximation}
Herbert Robbins and Sutton Monro.
\newblock A Stochastic Approximation Method.
\newblock {\em The Annals of Mathematical Statistics}, 22(3):400--407, 1951.

\bibitem{Sny-particle}
C.~Snyder, T.~Bengtsson, P.~Bickel, and J.~Anderson.
\newblock Obstacles to high-dimensional particle filtering.
\newblock {\em Mon. Wea. Rev.}, 136:4629--4640, 2008.

\bibitem{ExKF}
{T. Song} and J.~{Speyer}.
\newblock A stochastic analysis of a modified gain extended kalman filter with
  applications to estimation with bearings only measurements.
\newblock {\em IEEE Transactions on Automatic Control}, 30(10):940--949, 1985.

\bibitem{Tong_EnKF}
X. Tong, A.J. Majda, and D. Kelly.
\newblock Nonlinear stability and ergodicity of ensemble based {K}alman
  filters.
\newblock {\em Nonlinearity}, 29(2):657--691, 2016.

\bibitem{vanLeeuwen}
P.~J. van Leeuwen.
\newblock Nonlinear data assimilation in geosciences: an extremely efficient
  particle filter.
\newblock {\em Q. J. Roy. Meteor. Soc.}, 136(653):1991--1999, 2010.

\bibitem{zakai}
M. Zakai.
\newblock On the optimal filtering of diffusion processes.
\newblock {\em Z. Wahrscheinlichkeitstheorie und Verw. Gebiete}, 11:230--243,
  1969.

\bibitem{Zhang_Zakai}
H.~Zhang and D.~Laneuville.
\newblock Grid based solution of zakai equation with adaptive local refinement
  for bearing-only tracking.
\newblock {\em IEEE Aerospace Conference}, 2008.

\end{thebibliography}

\end{document}